\documentclass{amsart}

\usepackage{amssymb}
\usepackage{eepic}
\usepackage{color}
\usepackage{mathtools,xparse}
\usepackage[linktocpage=true]{hyperref}
\usepackage[utf8x]{inputenc}
\allowdisplaybreaks

\usepackage{comment}
\usepackage[usenames,dvipsnames,table]{xcolor}
\usepackage{tikz}

\definecolor{cranegreen}{cmyk}{1,0.76,0.94,0}

\newtheorem{thm}{Theorem}[section]
\newtheorem{pro}[thm]{Problem}
\newtheorem{cor}[thm]{Corollary}

\newtheorem{rem}[thm]{Remark}

\theoremstyle{definition}

\newtheorem{conj}[thm]{Conjecture}

\numberwithin{equation}{section}

\newcommand{\Z}{\mathbf{Z}}

\newcommand{\N}{\mathbf{N}}

\newcommand{\R}{\mathbf{R}}

\newcommand{\C}{\mathbf{C}}
\newcommand{\SL}{\mathrm{SL}}

\DeclareMathOperator{\Aff}{Aff}

\DeclareMathOperator{\supp}{supp}

\newcommand{\T}{\mathbf{T}}
\newcommand{\cL}{\mathcal{L}}
\newcommand{\B}{\mathcal{B}}

\def\ten{\otimes}
\def\G{\mathrm{G}}
\def\1{\mathbf{1}}
\def\H{\mathcal{H}}

\def\RR{\mathcal{R}}

\newcommand{\fin}{\hspace*{\fill} $\square$ \vskip0.2cm}

\renewcommand{\thefootnote}{\dag}

\begin{document}
\null

\vskip5pt

\null

\begin{center}
{\huge The impact of Schur multipliers in \\ harmonic analysis and operator algebras}

\vskip10pt

{\sc {Javier Parcet}}
\end{center}

\title[Fourier and Schur multipliers]{}

\author[Javier Parcet]{}

\maketitle

\null

\vskip-48pt

\null

\begin{center}
{\large {\bf Abstract}}
\end{center}

\vskip-28pt

\null

\begin{abstract}
Schur multipliers are basic linear maps on matrix algebras. Their close albeit still intriguing connection with Fourier multipliers establishes a powerful bridge between harmonic analysis and operator algebras. In this paper, we survey their growing impact over the past 15 years. Particular attention will be drawn to recent bounds on Schatten $p$-classes, with far-reaching applications in harmonic analysis on group von Neumann algebras and operator rigidity phenomena for higher-rank Lie groups and lattices. Key novelties arise from new insights into nonToeplitz Schur multipliers and unprecedented connections with highly singular operators from Euclidean harmonic analysis. 
\end{abstract}

\maketitle

\addtolength{\parskip}{+1ex}

\section*{\bf \large Introduction}

Inspired by Heisenberg's matrix mechanics and connections with ergodic theory and dynamical systems, the theory of von Neumann algebras emerged in the 1930s as a noncommutative form of measure theory. Here functions are replaced by linear operators on Hilbert spaces, which lack a commutative product. Harmonic analysis over von Neumann algebras extends beyond traditional \lq\lq noncommutative  harmonic analysis\rq\rq${}$ which prioritizes the (commutative) measure spaces formed by a nonabelian group and its Haar measure. While noncommutative $L_p$-spaces have been widely studied from a functional analytic viewpoint, they have certainly been underexploited in harmonic analysis, since it demands a genuine interaction across several fields. A key turning point in the late 20th century ---notably under Pisier's influence--- led to new noncommutative $L_p$-techniques, thanks to the development of operator space theory as well as quantum probability. Sums of independent or free random variables, square and maximal functions, martingale inequalities, singular integral operators and more were investigated thereafter.


A strong motivation for this form of harmonic analysis comes from structural properties of group von Neumann algebras, with major implications in geometric group theory and the classification of nonamenable factors. Haagerup's pioneering work on free groups and rank-one lattices \cite{DCH,CH,H} encoded deep geometric aspects of these groups in terms of Fourier approximation properties. In 2011 Lafforgue and de la Salle's theorem \cite{LdlS} sparked a powerful revival of Haagerup's methods, with a remarkable analysis of the failure of $L_p$-approximations in higher ranks and a smart use of Schur multipliers. These are linear maps on matrix algebras |first investigated by Issai Schur \cite{Schur} in 1911| whose definition is rather simple for finite matrices $$S_M(A) = \Big( M(j,k)A_{j,k} \Big)_{1 \le j,k \le N} \quad \mbox{for} \quad M \colon \{1,2,\ldots,N\} \times \{1,2,\ldots,N\} \to \C.$$ 

In this paper, we survey a tight relation between Fourier and Schur multipliers along with some new fundamental inequalities for them, and a variety of far-reaching applications in matrix and group von Neumann algebras. The extra flexibility to manipulate Schur multipliers compared to their Fourier peers is a crucial aspect of these achievements. Finally, we briefly review some implications in the context of Lafforgue/de la Salle's theorem and how new connections with Euclidean harmonic analysis could be relevant towards Connes' rigidity conjecture. Along the way, we shall also explore related topics and open problems. Overall, we aim to expose why Schur multipliers have become key tools in addressing different challenges in harmonic analysis, operator algebras and geometric group theory.    

\section{\bf \large New inequalities for Schur multipliers}\label{sec:schur}

Let $(\Omega,\mu)$ be a $\sigma$-finite measure space. Given $p \ge 1$, the Schatten class $S_p(\Omega)$ is the Banach space of all bounded linear maps $A \colon L_2(\Omega) \to L_2(\Omega)$ with finite norm $\|A\|_{S_p} = (\mathrm{tr} |A|^p)^{1/p}$, where $|A|^p$ arises from $A$ by functional calculus. The Schatten class $S_2(\Omega)$ is the space of Hilbert-Schmidt operators on $L_2(\Omega)$ and coincides with $L_2(\Omega \times \Omega)$ by identifying  $A \in S_2(\Omega)$ with its kernel $K_A$
\[Af(x) = \int_\Omega K_A(x,y) f(y) d\mu(y).\] 
Given $M \colon \Omega \times \Omega \to \C$ measurable, its \emph{Schur $S_p$-multiplier} is defined (when it exists) as the unique bounded  operator $S_M$ on $S_p(\Omega)$ assigning $A \in S_2(\Omega) \cap S_p(\Omega)$ to the operator $S_M(A)$ with kernel $M(x,y) K_A(x,y)$. In what follows, we shall formally identify $K_A(x,y)$ with a matrix $(A_{x,y})$. Next, $S_M$ is completely $S_p$-bounded when additionally $S_M \otimes \mathrm{id}_{S_p(\Gamma)}$ extends to a bounded map on $S_p(\Omega \times \Gamma)$ for any countable index set $\Gamma$. Given $1 < p \neq 2 < \infty$, Pisier conjectured in 1998 the existence of Schur $S_p$-multipliers failing complete boundedness, but no examples are known. We refer to \cite[Section 1]{LdlS} to see why no such examples exist  when $\Omega$ has no atoms and for a rather complete presentation of basic properties of Schur multipliers. 

\subsection{Fourier-Schur transference}

Fourier multiplier theory is a cornerstone in harmonic analysis. Given $m \colon \Z \to \C$ bounded, the Fourier multiplier $T_m$ is densely defined on square-integrable functions $f \colon \T \to \C$ by pointwise multiplication on their Fourier coefficients $\widehat{T_m} f(k) = m(k) \widehat{f}(k)$. It is \emph{completely $L_p$-bounded} when $T_m \otimes \mathrm{id}_{S_p(\Gamma)}$ extends to a bounded map on $L_p(\T; S_p(\Gamma))$, a space of matrix-valued functions. We set $$\big\| T_m \colon L_p(\T) \to L_p(\T) \big\|_{\mathrm{cb}} := \sup_{\Gamma} \big\| T_m \otimes \mathrm{id}_{S_p(\Gamma)} \colon L_p(\T; S_p(\Gamma)) \to L_p(\T; S_p(\Gamma)) \big\|.$$ A similar definition applies for symbols $m \colon \mathrm{G} \to \C$ on any locally compact abelian group $\mathrm{G}$ and Fourier multipliers on $L_2$-functions over its dual group. Complete boundedness is a strengthening in the category of operator spaces of Banach space boundedness \cite{PisAst,P2}. A great portion of Euclidean Fourier $L_p$-multiplier theory holds verbatim in the cb-setting. The Marcinkiewicz, H\"ormander-Mikhlin and Carleson-Sj\"olin multiplier theorems or de Leeuw's transference theorems are some illustrations. The cb-validity of the Littlewood-Paley-Rubio de Francia theorem \cite{RdFLP} in its full form is a notable open problem. The main known Fourier $L_p$-multipliers which fail cb-boundedness are due to Pisier and arise from $\Lambda_p$-sets whose Hankel extensions fail unconditionality in the matrix unit system \cite[Chapter 8]{PisAst}.    

Pisier's counterexample above relies on a tight relation between the trigonometric system of characters $\{\chi_g \colon \widehat{\G} \to \T \ \mathrm{s.t.} \ g \in \mathrm{G}\}$ and the system of matrix units $\{e_{g,h} \colon g,h \in \mathrm{G}\}$. To be more precise, the $*$-homomorphism $$\int_\mathrm{G} \widehat{f}(g) \chi_g(\cdot) \, d\mu(g) \longmapsto \Big( \widehat{f}(gh^{-1}) \Big)$$ is well-defined for $\widehat{f} \in \mathcal{C}_c(\G)$ and extends to a map $\Lambda \colon L_\infty(\widehat{\mathrm{G}}) \to \B(L_2(\G))$ by weak-$*$ density. Here, $\mu$ denotes the Haar measure of $\mathrm{G}$ and this suggests a relation between Fourier multipliers in the dual group $(\widehat{\mathrm{G}},\widehat{\mu})$ and certain class of Schur multipliers over $(\mathrm{G},\mu)$. This is indeed confirmed by the following result. 

\begin{thm} \emph{(Fourier-Schur transference I)}\textbf{.} \label{ThmFS1} Let $1 \le p \le \infty$ and consider a locally compact abelian group $\mathrm{G}$. Assume that $m \colon \mathrm{G} \to \C$ defines a completely $L_p$-bounded Fourier multiplier and set $M(g,h) = m(gh^{-1})$. Then $$\big\| S_M \colon S_p(\G) \to S_p(\G) \big\|_{\mathrm{cb}} = \big\| T_m \colon L_p(\widehat{\mathrm{G}}) \to L_p(\widehat{\mathrm{G}}) \big\|_{\mathrm{cb}}.$$
\end{thm}

Schur multipliers defined with symbols of the form $M(g,h) = m(gh^{-1})$ in some topological group are known as \emph{Herz-Schur} or \emph{Toeplitz} multipliers. Theorem \ref{ThmFS1} for $p=\infty$ follows from the elementary identity $\Lambda \circ T_m = S_M \circ \Lambda$ and was originally noted by Bo\.zejko/Fendler \cite{BF0} in 1984. The $L_p$-case requires to suitably modify the map $\Lambda$, which relies in turn on certain limiting/averaging process. It holds for nonabelian amenable groups as well (see below). The first result in this direction is due to Neuwirth/Ricard \cite{NR}, who proved it for amenable discrete groups. The general case was finally settled by Caspers/de la Salle in \cite{CS}. This important transference theorem |along with its nonabelian form in Theorem \ref{ThmFS2}| allows to encode Fourier multipliers as Schur multipliers with Toeplitz symbols, which will be particularly useful in the rest of this paper. 

NonToeplitz multipliers are no longer tied to Fourier multipliers and Theorem \ref{ThmFS1} shows that any result valid for arbitrary Schur multipliers can be understood as a \emph{nontrigonometric extension of a Fourier multiplier theorem}|the one which follows by restriction to Toeplitz symbols. The great flexibility to cut, restrict or deform Schur multipliers |\cite[Section 1]{LdlS} and \cite[Lemma 2.1]{PST}| makes nonToeplitz analysis extremely versatile, since Fourier multipliers (equivalently the Toeplitz subclass of Schur multipliers) are far more rigid. NonToeplitz multipliers are much less understood though and their $S_p$-boundedness is certainly mysterious. The Grothendieck celebrated inequality is closely connected to a characterization of the operator boundedness of Schur multipliers \cite{Gr,PisBAMS}. NonToeplitz multipliers also played a key role to solve the longstanding Krein's problem on operator-Lipschitz functions \cite{PS}. We will return to these results below. Other results and applications can be found in \cite{AK,Bennett,Ha,JP2,PisS}. In spite of these and other results in literature, sufficient conditions for $S_p$-boundedness were rather limited so far.   

\subsection{H\"ormander-Mikhlin-Schur multipliers} 

Let us momentarily fix $(\Omega,\mu)$ as the Euclidean space $\R^n$ equipped with its Lebesgue measure. By Fourier-Schur transference and since the H\"ormander-Mikhlin theorem \cite{Ho,Mi} holds as well in the cb-setting, we may rewrite it as follows for $M(x,y) = m(x-y)$
\begin{equation} \tag{HM} \label{Eq-HM}
\hskip-11pt 
\big\| S_M \colon S_p(\R^n) \to S_p(\R^n) \big\|_{\mathrm{cb}} \, \lesssim \, \frac{p^2}{p-1} \hskip-2pt \sum_{|\gamma| \le [\frac{n}{2}] +1} \Big\| |\xi|^{|\gamma|} \partial_\xi^\gamma m(\xi)  \Big\|_\infty.
\end{equation}
The singularity at $0$ for $\gamma \neq 0$ mirrors fundamental singular integrals in harmonic analysis which also appear in PDEs, differential geometry or fluid mechanics. The order $[\frac{n}{2}]+1$ is optimal. Radial $L_p$-multipliers satisfy reciprocally $|\xi|^{|\gamma|} \partial_\xi^\gamma m(\xi) \in L_\infty$ up to order $[\frac{n-1}{2}]$ for arbitrarily large (finite) $p$. Characterizing Fourier multipliers in $L_p$ for $p \neq 1,2,\infty$ is simply out of reach and the H\"ormander-Mikhlin theorem is one of the finest results for Fourier $L_p$-multipliers and $1<p<\infty$.   

In 2019, Mikael de la Salle conjectured a nontrigonometric/noncommutative form of the H\"ormander-Mikhlin theorem asking for a regularity condition on nonToeplitz symbols $M \colon \R^n \times \R^n \to \C$ outside the diagonal (admitting certain singularity on it) which implies the complete
$S_p$-boundedness of the associated Schur multiplier $S_M$ for $1 < p < \infty$. The following addresses de la Salle's question. 

\begin{thm} \label{ThmHMS} \emph{(H\"ormander-Mikhlin-Schur multipliers)}\textbf{.}
Let $1 < p < \infty$ and let $M \in \mathcal{C}^{[\frac{n}{2}]+1}(\R^{2n} \setminus \{x=y\})$ be a smooth symbol outside the diagonal. Then, we get
$$\big\| S_M \hskip-3pt: S_p(\R^n) \hskip-1pt \to \hskip-1pt S_p(\R^n) \big\|_{\mathrm{cb}} \hskip-1pt \le C_p \hskip-6pt \sum_{|\gamma| \le [\frac{n}{2}] +1} \hskip-3pt \Big\| |x-y|^{|\gamma|} \Big\{ \big| \partial_x^\gamma M(x,y) \big| + \big| \partial_y^\gamma M(x,y) \big| \Big\} \Big\|_\infty.$$ \vskip-10pt \noindent This recovers \eqref{Eq-HM} for Toeplitz symbols and the constant $C_p$ still behaves like $\frac{p^2}{p-1}$.
\end{thm} 

This nonToeplitz extension of the H\"ormander-Mikhlin theorem was established in \cite{CGPT1} and gives a rather easy-to-check criterion for cb-boundedness on Schatten $p$-classes with multiple applications outlined below. 

\noindent \textbf{Sketch of the proof.} 
Let $\RR = L_\infty(\R^n) \bar\ten \mathcal{B}(L_2(\R^n))$ be the von Neumann algebra of matrix-valued Euclidean functions and set $uf(x,z) = \exp \big( 2\pi i \langle x,z \rangle \big) f(x,z)$, a unitary on $L_2(\R^n \times \R^n)$. Define the representation $$\pi:  \B(L_2(\R^n)) \ni A \mapsto u \big( \1 \otimes A \big) u^* \in L_\infty(\RR).$$ Then, this map is a $*$-homomorphism satisfying that $$\pi(A) = \Big( \exp \big( 2\pi i \langle \, \cdot \, , x-y \rangle \big) A_{x,y} \Big)$$ for $A$ in the weak-$*$ dense subspace $S_2(\R^n)$. Next, using Mei's operator-valued BMO spaces \cite{Mei07} we define a BMO space for matrix algebras as the weak-$*$ closure of $\pi(\B(L_2(\R^n)))$ in Mei's space $\mathrm{BMO}_\RR$. By a simple duality argument, the proof is then reduced to the following results of independent interest:
\vskip3pt
i) \textbf{An interpolation theorem} $$\big[ \mathrm{BMO}, S_2(\R^n) \big]_{\frac{2}{p}} \simeq_{\mathrm{cb}} S_p(\R^n) \quad \mbox{for} \quad 2 \le p < \infty.$$ We omit here the details of this proof. The equivalence constant $c_p \approx p$ as $p \to \infty$.

\vskip2pt

ii) \textbf{Transference to twisted multipliers.} Noncommutative BMO spaces come as intersection of row and column forms $\mathrm{BMO} = \mathrm{BMO}_c \cap \mathrm{BMO}_r$, see \cite{Mei07,PX}. The key novelty here is to decouple the endpoint inequality $S_\infty \to \mathrm{BMO}$ via two different  transferences of matrix inequalities into operator-valued ones. To do so, we consider 
$$\widetilde{T}_{M_r}(f) = \Big( T_{M_r(\cdot,y)}(f_{x,y}) \Big) \quad \mbox{and} \quad \widetilde{T}_{M_c}(f) = \Big( T_{M_c(x, \cdot)}(f_{x,y}) \Big),$$ with $M_r(x,y) = M(y-x,y)$ and $M_c(x,y) = M(x,x-y)$. $\widetilde{T}_{M_r}$ and $\widetilde{T}_{M_c}$ are \lq\lq twisted Fourier multipliers\rq\rq${}$ acting on matrix-valued functions entrywise by Fourier multipliers which change with the matrix entry. Then $\pi \circ S_M = \widetilde{T}_{M_c} \circ \pi$ and we get

\begin{eqnarray*} 
\big\| S_M \colon S_\infty(\R^n) \to \mathrm{BMO} \big\| \!\! & \approx & \!\! \sum_{\dag \in \{r,c\}} \big\| S_M \colon S_\infty(\R^n) \to \mathrm{BMO}_\dag \big\| 
\\ [-1pt] \!\! & \le & \!\! \sum_{\dag \in \{r,c\}} \big\| \widetilde{T}_{M_\dag} \colon L_\infty(\RR) \to \mathrm{BMO}_{\RR}^\dag \big\| 
\\ [-2pt] \!\! & \le & \!\! \hskip-1pt \sum_{\star \in \{x,y\}} \sum_{|\gamma| \le [\frac{n}{2}] +1} \Big\| |x-y|^{|\gamma|} \partial_\star^\gamma M(x,y) \Big\|_\infty. 
\end{eqnarray*} 

\vskip2pt

iii) \textbf{Noncommutative Calder\'on-Zygmund methods.} 
Theorem \ref{ThmHMS} follows from points i) - ii). In order to justify the last inequality above, we should identify the kernel of twisted Fourier multipliers. We formally have
$$\widetilde{T}_{M_c}(f)(z) = \int_{\R^n} K_c(z-w) \cdot f(w) \, dw$$
with $K_c: \R^n \setminus \{0\} \to \mathcal{B}(L_2(\R^n))$ the \lq\lq diagonal-valued" function $$(K_c(z)\varphi)(x) = \big[ M_c(x,\cdot) \big]^{\vee}(z) \varphi(x).$$ A careful analysis of these kernels shows that $K_c$ should be understood as an operator-valued distribution which agrees with a locally integrable operator-valued function on $\R^n \setminus \{0\}$. In particular, the representation above is meaningful when $z \notin \mathrm{supp}_{\R^n} f$, the Euclidean support of the matrix-valued function $f$. Under this assumption (standard in CZ theory), the above inequalities follow after refining previous methods from noncommutative CZ theory \cite{CCP,JMP1,Pa1}, we omit details. \fin

\begin{rem} \label{RemNCCZ}
\emph{Noncommutative Calderón-Zygmund theory has been quite useful and influential over the past 15 years. The lack of a weak type $(1,1)$ inequality for Calder\'on-Zygmund operators acting on matrix-valued functions was noticed by Pisier and Xu circa 2005. This was addressed in the seminal paper \cite{Pa1}, where noncommutative martingale theory was used to produce a noncommutative form of Calder\'on-Zygmund decomposition. CZ operators in crossed product von Neumann algebras (Euclidean measure spaces) gave rise in \cite{JMP1} to the first H\"ormander-Mikhlin type theorems in group von Neumann algebras. Purely noncommutative scenarios avoiding Euclidean spaces in tensor or crossed products were also investigated. In \cite{GJP} the archetype manifolds in noncommutative geometry were considered. This includes noncommutative tori and the Heisenberg-Weyl algebra, along with other quantum Euclidean spaces which appear in quantum field theory, string theory or quantum probability. A general CZ theory for von Neumann algebras was developed in \cite{JMPX} under algebraic assumptions in terms of Markov processes which replace standard metric assumptions. A simpler CZ decomposition including nondoubling measures was found in \cite{CCP}. Noncommutative CZ decompositions have been crucial in the solution of Nazarov-Peller's conjecture \cite{CPSZ} and Cadilhac/Wang's remarkable ergodic theorem \cite{CW}, extending \cite{HLW,JX}. Further developments and applications of noncommutative Calder\'on-Zygmund theory appear in \cite{Ca,CR,HLX,HLM,PRS}. Still two open problems have remained open. One is whether H\"ormander's kernel condition implies the weak type $L_1$ inequality for matrix-valued CZ operators, see \cite{CCP} for valid slightly stronger assumptions. The other is the weak type $L_1$ extension of Theorem \ref{ThmHMS}, with potential applications in group von Neumann algebras. Also commutator estimates are missing. New insights for them could improve the Calder\'on-Zygmund decomposition for diagonal-valued kernels.}
\end{rem}

\null

\vskip-25pt

\null

\begin{center}
\includegraphics[scale=0.45]{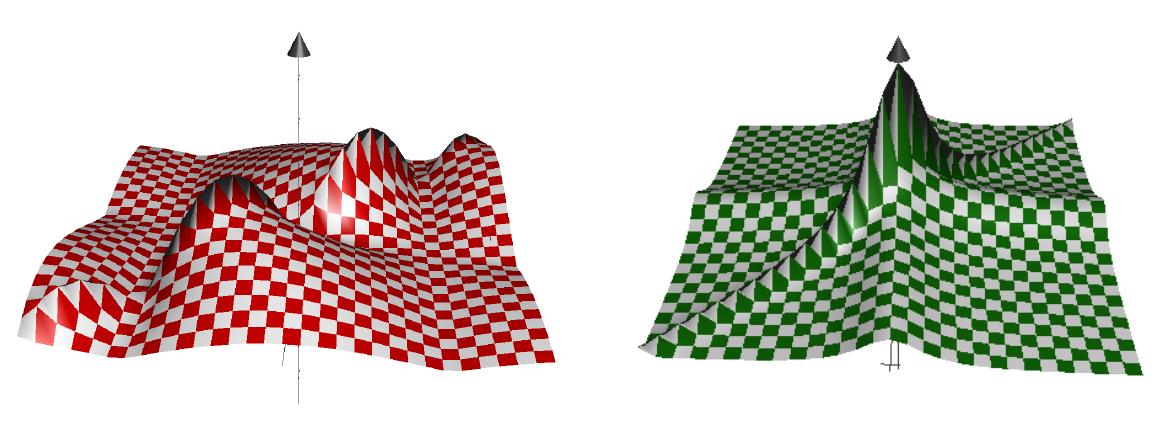}
\end{center}

\vskip-30pt

\null

\begin{center}
{\scriptsize \textbf{Figure 1.} Examples of NonToeplitz H\"ormander-Mikhlin-Schur multipliers in $\R \times \R$} \\ [-3pt]
{\scriptsize Any Toeplitz symbol would be forced to be constant at $x=y+\alpha$ for all $\alpha \in \R$, unlike above}
\end{center}

\vskip5pt

We conclude our analysis of H\"ormander-Mikhlin-Schur multipliers with a brief overview of interesting related topics and applications in matrix algebras. Other implications in group von Neumann algebras will be explored in Section \ref{sec:fourier}.

A. \textbf{Marcinkiewicz-Schur multipliers.} Another landmark result in Fourier multiplier theory is Marcinkiewicz's theorem \cite{Gra}. Given $1 < p < \infty$, it ensures $L_p$-boundedness for bounded symbols $m \colon \R \to \C$ which have bounded variation over dyadic intervals. Their cb-boundedness was confirmed by Bourgain \cite{Bou} and the nonToeplitz extension was established in \cite{CLM} by Chuah, Liu and Mei. The argument shares ideas with Theorem \ref{ThmHMS} but was found independently. Given $M \colon \Z \times \Z \to \C$ bounded and setting $\mathcal{J}_m = \{k \in \Z \colon 2^{m-1} \le |k| < 2^m\}$, their result shows that the Schur multiplier $S_M$ is completely $S_p(\Z)$-bounded for $1 < p < \infty$ as long as the quantity below is finite $$\sup_{\begin{subarray}{c} j \in \Z \\ m \in \N \end{subarray}} \, \sum_{k \in \mathcal{J}_m} \big| M(j+k+1,j) - M(j+k,j) \big| + \big| M(j,j+k+1) - M(j,j+k) \big|.$$ Continuous versions in the real line $\R$ and higher-dimensional analogs are also given. A particularly nice remark there is that both terms above |as well as both terms $\partial_x^\gamma$ and $\partial_y^\gamma$ in Theorem \ref{ThmHMS}| are necessary, not just an artifact of the proof.

B. \textbf{Refining Arazy's conjecture.} Given $1 < p < \infty$, a Lipschitz function $f: \R \to \C$ and $A, B$ self-adjoint operators with $A-B \in S_p(\R)$, Potapov and Sukochev proved in \cite{PS} that $\|f(A) - f(B)\|_{S_p(\R)} \le C_p \|f\|_{\mathrm{Lip}} \|A-B\|_{S_p(\R)}$ for some constant $C_p$. This solves a longstanding problem posed by Krein in 1964. Their proof consisted in showing the validity of Arazy's 1982 stronger conjecture for divided differences $$M_f(x,y) = \frac{f(x)-f(y)}{x-y} \ \leadsto \ \big\| S_{M_f} \colon S_p(\R) \to S_p(\R) \big\| \le C_p \|f\|_{\mathrm{Lip}}.$$ This statement is inherently nonToeplitz since its Toeplitz form imposes $f$ to be linear, which makes $M_f$ constant. Nazarov-Peller's conjecture alluded in Remark \ref{RemNCCZ} is the weak type $(1,1)$ analog of Krein's conjecture. H\"ormander-Mikhlin-Schur multipliers go far beyond Arazy's conjecture. Namely, in first place Theorem \ref{ThmHMS} provides a one-line proof of the main result in \cite{PS} with optimal constants. 
Also an improved version of Theorem \ref{ThmHMS} near $L_2$ á la Calder\'on-Torchinsky allows to generalize Arazy's conjecture to $\alpha$-H\"older divided differences for $0 < \alpha < 1$ |the case $\alpha = 1$ corresponding to Arazy's conjecture| further details in \cite{CGPT1}. Another application of H\"ormander-Mikhlin-Schur multipliers (Theorem \ref{ThmHMS}) was recently given in \cite{CaRe} for second-order divided differences. 

\begin{rem} \label{Rem-Riesz-Schur} \emph{A different proof of Theorem \ref{ThmHMS} was recently found in \cite{GPPR}. This combines a powerful technique from \cite{JMP2} |\emph{H\"ormander-Mikhlin multipliers appear as Littlewood-Paley averages of Riesz transforms over fractional laplacians}, which was somehow hidden in Bourgain's work \cite{Bou} as later noted by Ricard| with a nonToeplitz form of dimension-free estimates for Riesz transforms, which we shall comment below. Notably, the proof avoids CZ and probabilistic methods.}
\end{rem}

\begin{rem}
\emph{Extensions of Theorem \ref{ThmHMS} over locally compact groups appear in \cite[Theorem 4.1, Corollary 4.5]{CGPT2}, see also \cite[Remark 4.2]{CGPT2} for other measures spaces.}
\end{rem}

\subsection{The local geometry of idempotent Schur multipliers} 

Rigidity aspects of high rank lattices from \cite{LdlS} strongly motivated Theorem \ref{ThmHMS}, but there is still much to learn about less regular multipliers. A key point in \cite{dLdlS,LdlS} was a careful analysis of Schur multipliers over the $n$-sphere for symbols $M_\varphi(x,y) = \varphi(\langle x,y \rangle)$ depending on the inner product of its entries. More precisely, the boundedness of $S_{M_\varphi}$ on the Schatten class $S_p(\mathbf{S}^n)$ for $p > 2 + \frac{2}{n-1}$ imposes some H\"older regularity on $\varphi$. Can we admit less regular multipliers closer to $L_2$? How less? How close? The $S_p$-mapping properties of the \emph{spherical Hilbert transform} 
$$H_\mathbf{S}: A \mapsto \Big( -i \, \mathrm{sgn} \langle x,y \rangle A_{x,y} \Big)_{x,y}$$
is a very basic problem in this regard concerning jump discontinuities. Is $H_\mathbf{S}$ an $S_p$-bounded map for some $\frac{2n}{n+1} < p \neq 2 < \frac{2n}{n-1}$? Equivalently, we may consider the Schur multiplier $(1+iH_\mathbf{S})/2$ with symbol $\chi_\Sigma$ for $\Sigma = \{\langle x,y\rangle>0\}$. We could even consider other idempotent multipliers|whose symbols are characteristic functions of smooth domains. The analogy with Fefferman's celebrated theorem for the ball \cite{Fe} |before which unboundedness was only known for $p$ outside this range| is worth noting. Theorem \ref{ThmIdemp} solves this problem with a vast generalization of \cite{Fe}. 

Let $M$ be a differentiable manifold with the Lebesgue measure coming from any Riemmanian structure on it. Consider a $\mathcal{C}^1$-domain $\Sigma \subset M \times M$ so that its boundary $\partial \Sigma$ is a smooth hypersurface, which is locally represented by level sets of some real-valued $\mathcal{C}^1$-functions with nonvanishing gradients. We say that $\partial \Sigma$ is \emph{transverse} at a point $(x,y)$ when the tangent space of $\partial \Sigma$ at $(x,y)$ maps surjectively on each factor $T_x M$ and $T_y M$. In that case, both sections $$\partial \Sigma_x = \big\{ y' \in M \mid (x,y') \in \partial \Sigma \big\} \quad \mbox{and} \quad \partial \Sigma^y = \big\{ x' \in M \mid (x',y) \in \partial \Sigma \big\}$$ become codimension $1$ manifolds on some neighbourhood of $y$ and $x$ respectively.

\begin{thm} \label{ThmIdemp} \emph{(Geometry of idempotent Schur multipliers)}\textbf{.}
Let $1 < p \neq 2 < \infty$ and consider a $\mathcal{C}^1$-domain $\Sigma \subset M \times M$. \hskip-1pt Then the following statements are equivalent for any transverse point $(x_0,y_0) \in \partial \Sigma \hskip-2pt :$
\begin{enumerate}
\item\label{item:Spbounded}  \emph{$S_p$-boundedness}. The Schur multiplier $S_\Sigma$ with symbol $\chi_\Sigma$ is locally bounded on $S_p(M)$ around $(x_0,y_0)$. That is, $S_{\Sigma \cap (U \times V)}$ is $S_p(M)$-bounded for some neighbourhoods $U,V$ of $x_0, y_0$ in $M$.

\vskip2pt

\item\label{item:SteinCurvature} \emph{Zero-curvature condition}. There are neighbourhoods $U,V$ of $x_0, y_0$ in $M$ such that the tangent spaces $T_y (\partial \Sigma_{x_1})$ and $T_y (\partial \Sigma_{x_2})$ coincide for any pair of points $(x_1,y), (x_2,y) \in \partial \Sigma \cap (U \times V)$. 

\vskip2pt

\item\label{item:HilbertTransform} \emph{Triangular truncation representation}. There are neighbourhoods $U,V$ of the points $x_0, y_0$ in $M$ and $\mathcal{C}^1$-functions $f_1 \colon U \to \R$ and $f_2 \colon V \to \R$, such that we have $\Sigma \cap (U \times V) = \big\{(x,y) \in U \times V \colon f_1(x)>f_2(y) \big\}$.
\end{enumerate}
\end{thm} 

Theorem \ref{ThmIdemp} was proved in \cite[Theorem A]{PST} and characterizes the local geometry of Schur $S_p$-idempotents, which gives a nontrigonometric extension of Fefferman's ball multiplier theorem. It is rather surprising that similar phenomena holds on general manifolds, which lack to admit a Fourier transform connection. A moment of thought shows that it is only interesting when $\dim M \ge 2$, in line with \cite{Fe}. Before sketching the proof, a few remarks are in order:

i) \textbf{On Fourier idempotents.} Fefferman's landmark theorem \cite{Fe} disproves the $L_p$-boundedness of the idempotent Fourier multiplier given by the Euclidean ball for $p \neq 2$ and dimensions higher than one. The same argument works for Fourier idempotents over smooth Euclidean domains $\Omega$ admitting a boundary point of nonvanishing curvature. By Fourier-Schur transference, Fefferman's theorem corresponds in Theorem \ref{ThmIdemp} to $(M,\Sigma) = \big( \R^n, \big\{ (x,y) \colon x-y \in \Omega \big\} \big)$ for a Euclidean $\mathcal{C}^1$-domain $\Omega$. Every boundary point is trivially transverse in this case and this partly explains why transversality did not appear so far. 

ii) \textbf{On the global behavior.} A first consequence of Theorem \ref{ThmIdemp} is that the given equivalent properties do not depend on $1 < p \neq 2 < \infty$. It is important to insist here that the characterization is local. If the global aspects are taken into account, the situation is different. In the setting of discrete index sets, we know from \cite{Ha} that there are $S_p$-bounded Schur idempotents with $p \in 2 \Z_+$ and which fail to be $S_q$-bounded for $q>p$. Also, other examples for continuous index sets where the local theorem above fails to be global are given in \cite[Appendix A]{CPPR} as explained in \cite[Remark 2.7]{PST}. These domains are necessarily not relatively compact.  

iii) \textbf{On the notion of curvature.} If $M=\R^n$ and $( \mathbf{n}_1(x_0,y_0),\mathbf{n}_2(x_0,y_0) ) \perp \partial \Sigma$ at $(x_0,y_0)$, transversality means that both $\mathbf{n}_1, \mathbf{n}_2 \neq 0$ and zero-curvature that $\mathbf{n}_2(x_1,y) \parallel \mathbf{n}_2(x_2,y)$. Smoother domains $\Sigma$ can be locally described as level sets of a $\mathcal{C}^2$-function $\Phi$. In that case, zero-curvature becomes $u_x^\mathrm{t} ( \partial_{x_j} \partial_{y_k} \Phi(x_0,y_0) ) u_y = 0$ for all $u_\dag \perp \nabla_\dag \Phi(x_0,y_0)$ with $\dag = x,y$. The absence of 2nd order noncrossed derivatives here ($\partial_{x_j} \partial_{x_k}$ or $\partial_{y_j} \partial_{y_k}$) is justified since one-variable sets $\Sigma_r = \{(x, y) \colon x \in \Omega\}$ and $\Sigma_c = \{(x, y) \colon y \in \Omega\}$ lead to $S_p$-contractions regardless the geometry of $\Omega$. This flexible notion of curvature is crucial here and will be for Lie groups below.    

\noindent \textbf{Sketch of the proof.} The implication \eqref{item:Spbounded}$\Rightarrow$\eqref{item:SteinCurvature} follows from a matrix form of the celebrated Meyer's lemma \cite[Lemma 1]{Fe}, claiming that $L_p$-boundedness of the ball multiplier would imply certain square-function $L_p$-inequalities for families of half-space multipliers $$\widehat{H_{u}f}(\xi) = \chi_{\langle \xi, u \rangle > 0} \widehat{f}(\xi).$$ Its matrix form \cite[Lemma 2.3]{PST} is stated in the terminology of Theorem \ref{ThmIdemp} as follows. Given $N \ge 1$, let $x_1, x_2, \ldots, x_N \in U$ and $y \in V$ such that $z_j = (x_j,y)$ are transverse points in $\partial \Sigma$. Define $u_j = \mathbf{n}_2(z_j)$ as we did above. Then 
\begin{equation} \tag{MS} \label{EqMatrixMeyer}
\Big\| \Big( \sum_{j=1}^N \big| H_{u_j}(f_j) \big|^2 \Big)^\frac12 \Big\|_{L_p(\R^n)} \le \big\| S_{\Sigma \cap (U \times V)} \big\|_{\B(S_p(\R^n))} \, \Big\| \Big( \sum_{j=1}^N |f_j|^2 \Big)^\frac12 \Big\|_{L_p(\R^n)}.
\end{equation}
Inequality \eqref{EqMatrixMeyer} is a new connection between Fourier and Schur multipliers which readily gives the implication \eqref{item:Spbounded}$\Rightarrow$\eqref{item:SteinCurvature}. Indeed, taking local charts we may assume that $M$ is $\R^n$. Then, by transversality the map $\mathbf{n}_2(z)/\|\mathbf{n}_2(z)\|$ is continuous around $(x_0, y_0)$ and the failure of \eqref{item:SteinCurvature} would give a continuous of distinct directions $u_j$ in the assumptions of \eqref{EqMatrixMeyer}. However, the resulting inequalities do not hold with constants independent of $N$ as it follows combining Fefferman's 2-dimensional Besicovitch type construction \cite{Fe} and de Leeuw's restriction theorem \cite{dL}. The proof of \eqref{EqMatrixMeyer} requires a more involved analytic argument. Indeed, the key point is the following identity for any transverse point $(x,y) \in \partial \Sigma$ and every $T \in G \hskip-1pt L_n(\R)$ satisfying $T^* \mathbf{n}_1(x,y) = - \mathbf{n}_2(x,y)$
$$\lim_{\varepsilon \to 0^+} \chi_\Sigma \big( x + \varepsilon T\xi, y + \varepsilon \eta \big) = \chi_{\langle \xi-\eta, \mathbf{n}_2(x,y) \rangle > 0} \quad \mbox{for a.e. $\xi,\eta \in \R^n$}.$$ This relates the geometry of $\Sigma$ with half-space multipliers. The argument is then completed with Fourier-Schur transference and the flexibility of Schur multipliers to be deformed preserving their $S_p$-mapping bounds \cite[Lemma 2.1]{PST}.  

The implication \eqref{item:SteinCurvature}$\Rightarrow$\eqref{item:HilbertTransform} is a purely geometric statement about hypersurfaces in product manifolds|Theorem 2.5 in \cite{PST}. It follows by noticing that both properties remain invariant under diffeomorphisms of product type $(x,y) \mapsto (\phi(x),\psi(y))$ and a nontrivial iterated use of the implicit function theorem. We failed to find a straightforward proof of this result and consulting with a few experts did not provide alternative ideas or references to consult. The implication \eqref{item:HilbertTransform}$\Rightarrow$\eqref{item:Spbounded} is easier. Again by well-known techniques summarized in \cite[Lemma 2.1]{PST}, one can locally deform the functions $f_1$ and $f_2$ with no effect in the $S_p$-mapping bounds of $S_\Sigma$. Thus, this implication follows from the classical $S_p$-boundedness of the triangular projection $(A_{j,k}) \mapsto (\chi_{j \geq k} A_{j,k})$ from \cite{Ma}. By Fourier-Schur transference, this can be related in turn to the $L_p$-boundedness of the Hilbert transform. \fin 

\null

\vskip-25pt

\null

\includegraphics[scale=0.75]{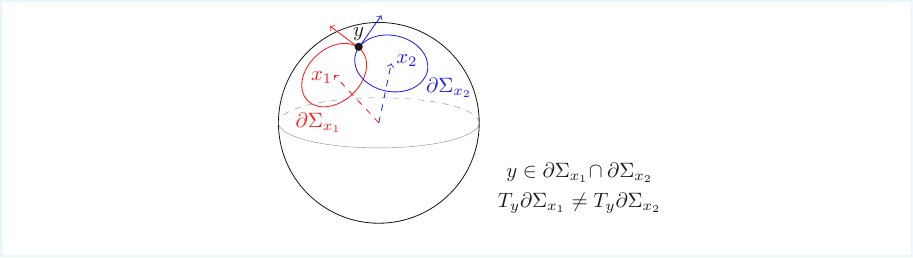}

\vskip-30pt

\null

\begin{center}
{\scriptsize \textbf{Figure 2.} Failure of \eqref{item:SteinCurvature} for spherical Hilbert transforms $H_{{\mathbf{S},\delta}}$} \\ [-3pt]
{\scriptsize Here $H_{{\mathbf{S},\delta}} = S_\Sigma$ with $\Sigma = \big\{(x,y) \in \mathbf{S}^n \times \mathbf{S}^n : \langle x,y \rangle > \delta \big\}$ for $n = 2$}
\end{center}

\vskip5pt

A few immediate outcomes also arise from Theorem \ref{ThmIdemp}. First, the zero-curvature condition \eqref{item:SteinCurvature} could be equivalently formulated using tangent spaces $T_x (\partial \Sigma^{y_1})$ and $T_x (\partial \Sigma^{y_2})$ of $y$-sections instead for any $(x,y_1), (x,y_2) \in \partial \Sigma \cap (U \times V)$. Also, the implication \eqref{item:Spbounded}$\Rightarrow$\eqref{item:HilbertTransform} shows that (up to $\mathcal{C}^1$-diffeomorphims of product type) the triangular projection is the only local model for Schur $S_p$-idempotents. Last, the transversality assumption is crucial in our proof for \eqref{item:Spbounded}$\Rightarrow$\eqref{item:SteinCurvature}$\Rightarrow$\eqref{item:HilbertTransform}, but it is still unclear whether it is necessary. Theorem \ref{ThmIdemp} trivially holds over the interior of the set of nontransverse points in the relative topology of $\partial \Sigma$, but it is likely that other transverse points in domains with nonanalytic boundary must be removed.

\begin{rem} 
\emph{Theorem \ref{ThmIdemp} also yields strong results on Lie groups, see Section \ref{sec:fourier}.}
\end{rem}

\begin{rem} \label{Rem-SphericalH} 
\emph{Spherical Hilbert transforms are $L_p$-unbounded in dimensions $n \ge 2$ for $p \neq 2$, as it follows from Theorem \ref{ThmIdemp} and illustrated in Figure 2. In particular spherical $S_p$-bounded multipliers with symbol $M_\varphi(x,y) = \varphi(\langle x,y \rangle)$ cannot include jump discontinuities for any $p \neq 2$. How regular must be $\varphi$ when $p$ approaches 2? This is a very subtle problem which puts together highly singular operators from Euclidean harmonic analysis with Connes' rigidity, see Section \ref{sec:rigidity}.}
\end{rem}

\section{\bf \large Harmonic analysis in group von Neumann algebras}\label{sec:fourier}

Let $\G$ be a locally compact topological group equipped with its left Haar measure $\mu$ and consider its left regular representation $\lambda \colon \G \to \mathcal{U}(L_2(\G))$, which is given by $\lambda_g (\varphi)(h) = \varphi(g^{-1}h)$. The group von Neumann algebra $\cL \G$ is the weak-$*$ closure in $\B(L_2(\G))$ of the space of Fourier expansions $$\Big\{ f = \int_\G \widehat{f}(g) \lambda_g \, d\mu(g)\colon \widehat{f} \in \mathcal{C}_c(\G) \Big\}.$$ Note that each such $f$ acts as a convolution operator $\varphi \mapsto \widehat{f} * \varphi$ on $L_2(\G)$. When $\G$ is abelian, $\cL \G$ may be easily identified with the $L_\infty$-space over the dual group of $\G$. The $L_p$-theory of $\cL \G$ is quite elementary for $\G$ unimodular. In that case, $\cL \G$ carries a natural trace $\tau \colon \cL \G  \to \C$ |known as Plancherel trace| determined by $$\tau(f^*f) = \int_\G \big| \widehat{f}(g) \big|^2 \, d\mu(g) \quad \mbox{for} \quad \widehat{f} \in L_2(\G).$$ If $e$ is the unit of $\G$, the standard identity $\tau(f) = \widehat{f}(e)$ holds for $\widehat{f} \in \mathcal{C}_c(\G) * \mathcal{C}_c(\G)$. Given $1 \le p < \infty$, the noncommutative $L_p$-space over the group algebra $\cL \G$ will be denoted $L_p(\cL \G)$ and is defined as the completion of $\{f \in \cL \G \colon \|f\|_p < \infty\}$ for the norm $\|f\|_p = \tau( |f|^p)^{1/p}$. Here, the Hilbert space operator $|f|^p = (f^*f)^{p/2}$ arises by functional calculus and $L_\infty(\cL \G)$ is just $\cL \G$ equipped with its operator norm, see e.g. \cite[Section 2]{CPPR} for further details. We refer to Pisier/Xu's excellent survey \cite{PX2} for an overview on the structure and properties of noncommutative $L_p$-spaces. A bounded measurable function $m \colon \G \to \C$ defines a Fourier $L_p$-multiplier when the map $$f \mapsto \int_\G m(g) \widehat{f}(g) \lambda_g \, d\mu(g)$$ extends to a bounded operator $T_m \colon L_p(\cL \G) \to L_p(\cL \G)$. Complete $L_p$-boundedness imposes additionally that $T_m \otimes \mathrm{id}_{S_p(\Gamma)}$ is bounded on the matrix amplification $L_p(\cL \G \bar\otimes \B(\ell_2(\Gamma))) = L_p(\cL \G;S_p(\Gamma))$ for any countable index set $\Gamma$. The definitions of $L_p(\cL \G)$ and Fourier $L_p$-multiplier above are more involved for nonunimodular groups. The reader may consult \cite[Section 2.3]{PST} and the references therein. 

\begin{rem}
\emph{Fourier multipliers on noncommutative $L_p$-spaces over group von Neumann algebras place the reference group in the frequency side. Motivated by problems in operator algebras, ergodic theory or geometric group theory, they are perhaps less known in harmonic analysis that other Fourier multipliers which act on classical $L_p$-spaces over type I topological groups, where frequencies appear as irreducible representations. Both (dual) settings coexists under a generalized form of Pontryagin duality in the general context of quantum groups \cite{KVPNAS}.}
\end{rem} 

\subsection{More on Fourier-Schur transference}

A locally compact group $\G$ is called amenable when it admits a left-invariant mean. In other words, when there exists a norm 1 positive preserving map $\Phi \colon L_\infty(\G) \to \R$ satisfying that $\Phi(\lambda_g (\varphi)) = \Phi(\varphi)$ for every $(g,\varphi) \in \G \times L_\infty(\G)$. Introduced by John von Neumann, amenable groups are characterized as those groups not admitting paradoxical decompositions á la Banach-Tarski. They admit many other equivalent definitions|one will be crucial in Section \ref{sec:rigidity} below. Using the above notion of Fourier multiplier, the result given below extends Theorem \ref{ThmFS1} to amenable groups and provides local variants for nonamenable ones. Both statements will be very useful in what follows. 

\begin{thm} \emph{(Fourier-Schur transference II)}\textbf{.} \label{ThmFS2} Let $1 \le p \le \infty$ and consider a locally compact group $\mathrm{G}$. Assume that $m \colon \mathrm{G} \to \C$ defines a completely $L_p$-bounded Fourier multiplier on its group algebra $\cL \G$ and set $M(g,h) = m(gh^{-1})$. Then, the following transference results hold$\hskip1pt :$ 
\begin{itemize}
\item[i)] If $\mathrm{G}$ is amenable, we have $$\big\| S_M \colon S_p(\G) \to S_p(\G) \big\|_{\mathrm{cb}} = \big\| T_m \colon L_p(\cL \G) \to L_p(\cL \G) \big\|_{\mathrm{cb}}.$$ Moreover, the upper inequality $\le$ holds for nonamenable groups as well.

\vskip3pt

\item[ii)] If $\mathrm{G}$ is nonamenable and $m \colon \mathrm{G} \to \C$ is compactly supported by $\Omega \subset \G$ $$\hskip39pt \big\| T_m \colon L_p(\cL \G) \to L_p(\cL \G) \big\|_{\mathrm{cb}} \le C_p(\G, \Omega) \big\| S_M \colon S_p(\G) \to S_p(\G) \big\|_{\mathrm{cb}}.$$ That is, Fourier-Schur transference holds locally for nonamenable groups. 
\end{itemize}
\end{thm}

Theorem \ref{ThmFS2} i) was proved for discrete amenable groups in \cite{NR} and generalized in \cite{CS}. Its local nonamenable form was first proved in \cite{PRS} for unimodular groups and $p \in 2 \Z_+$ and later generalized in \cite{PST}. As for abelian groups, this establishes a profound although still incomplete relation between the trigonometric and matrix unit systems associated to $\G$. Understanding it in further detail is in the root of several challenges at the interface of harmonic analysis and operator algebras.   

\noindent \textbf{Sketch of the proof.} The upper inequality in i) is rather simple. In case $\G$ is discrete, we just set $\mathbf{u} = \sum_{g \in \G} e_{g,g} \otimes \lambda_g$ and observe that unitary conjugation gives 
\begin{eqnarray*}
\big\| S_M(A) \big\|_{S_p(\G)} \!\!\! & = & \!\!\! \big\| \mathbf{u} \big( \1 \otimes S_M(A) \big) \mathbf{u}^* \big\|_{L_p(\cL \G; S_p(\G))} \\ \!\!\! & = & \!\!\! \big\| (T_m \otimes \mathrm{id}) \big( \mathbf{u} \big( \1 \otimes A \big) \mathbf{u}^*\big) \big\|_{L_p(\cL \G; S_p(\G))},
\end{eqnarray*}
from which the assertion follows after matrix amplification. The general case in \cite{CS} requires a careful reformulation using the subtle definition of Fourier multiplier for nonunimodular groups. The lower estimate in i) follows from ii), which gives $C_p(\G,\Omega) = 1$ for $\G$ amenable and $\Omega = \G$. The proof of ii) for $\G$ unimodular and $p \in 2\Z_+$ is sketched as follows. By translation invariance, we may assume that $\Omega$ is a relatively compact neighborhood of the identity. Then, there exists a constant $$0 \le \delta_\G(\Omega) := \inf_{\begin{subarray}{c} \|\phi\|_2 = 1 \\ \phi: \G \to \R_+ \end{subarray}} \sup_{g \in \Omega} \, \frac12 \int_\G \big| \phi(gh) - \phi(h) \big|^2 d\mu(h) < 1.$$ Using F\o lner sequences, it is easily seen that $\delta_\G(\Omega) = 0$ when $\G$ is amenable. In nonamenable groups, we get $\delta_\G(\Omega) \approx 0$ for $\Omega$ small enough and $\delta_\G(\Omega) \to 1$ as $\Omega \to \G$. This suggests that $\delta_\G(\Omega)$ measures the nonamenability of $\G$ relative to $\Omega$. Next, consider the (cb-contractive) maps $j_{p\phi}: L_p(\cL \G) \to S_p(\G)$ formally given by $$j_{p\phi}(f) = \Big( \phi(g)^{\frac{2}{p}} \widehat{f}(gh^{-1}) \Big)_{g,h \in \G}.$$ The key \cite[Lemma 1.3]{PRS} is to show for some $\phi \in \mathrm{B}_1(L_2^+(\G))$ that $$\supp \widehat{f} \subset \Omega \Longrightarrow \|f\|_p \le_{\mathrm{cb}} \frac{2}{1 - \delta_\G(\Omega_p)} \|j_{p\phi}(f)\|_{S_p(\G)} \quad \mbox{for} \quad \Omega_p = (\Omega \Omega^{-1})^{\frac{p}{2}}.$$ The assertion follows from it and $j_{p\phi}(T_m f) = S_M (j_{p\phi}f)$. The proof of ii) in the nonunimodular case \cite{PST} relies once more on involved definitions which we omit. \fin

\begin{rem}
\emph{The more efficient argument in \cite[Section 3]{PST} still fails to provide a constant $C_p(\G,\Omega)$ converging to $1$ as $p$ tends to $2$. This remains an open problem.}
\end{rem}

\begin{rem}
\emph{See \cite{CPPR} for noncommutative transference results á la de Leeuw \cite{dL}.}
\end{rem}

\subsection{Fourier multipliers with one-point singularities} 

The theory of Fourier $L_p$-multipliers on group von Neumann algebras has been intensively investigated over the past 15 years  \cite{C1,CGPT2,GJP0,JMP1,JMP2,JR,dLdlS,LdlS,MR,MRX,PRS,PRo,PST}. A great effort has been put in understanding noncommutative forms of  H\"ormander-Mikhlin criteria \cite{Ho,Mi}. We now revisit most of them using Schur multipliers from Theorem \ref{ThmHMS} and some variants. Related results will also be pointed along the way.  

\subsubsection{H\"ormander-Mikhlin criteria via cocycles} 

The classical H\"ormander-Mikhlin multiplier theorem establishes the $L_p$-boundedness of Euclidean Fourier multipliers under the regularity condition from \eqref{Eq-HM}, which allows one-point singularities in terms of the iterated derivatives of the symbol. Let now $\G$ be unimodular and equipped with its Haar measure $\mu$. The lack of a differential structure on $\G$ forces to find auxiliary ways to measure the regularity of a symbol $m \colon \G \to \C$. A broader interpretation of tangent spaces was exploited in \cite{JMP1} using cocycle maps, a standard tool from cohomology and representation theory. A $n$-dimensional cocycle of $\G$ is a map $\beta \colon \G \to \R^n$ together with an orthogonal action $\alpha \colon \G \to O_n(\R)$ satisfying the cocycle law $\alpha_g(\beta(h)) = \beta(gh) - \beta(g)$. The hope in \cite{JMP1} was that sufficient regularity conditions for a symbol $m \colon \G \to \C$ could exist in terms of the H\"ormander-Mikhlin condition for those Euclidean lifts $\widetilde{m}: \R^n \to \C$ satisfying $m = \widetilde{m} \circ \beta$. 

\null

\vskip-28pt

\null

\begin{center}
\hskip10pt \includegraphics[scale=0.8]{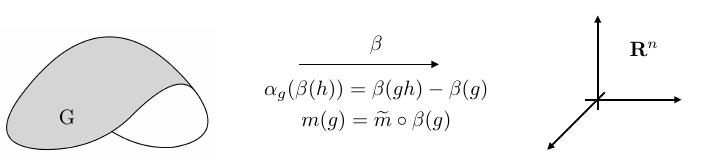}
\end{center}

\vskip-20pt

\null

\begin{center}
{\scriptsize \textbf{Figure 3.} Cocycle liftings of Fourier symbols $m: \G \to \C$ \\ [-3pt] Different cocycle maps provide a variety of sufficient conditions, see Remark \ref{Rem-Cocycles}}
\end{center}

\vskip5pt

The bimodularity of noncommutative Calderón-Zygmund methods required to lift $m$ through both \emph{left} and \emph{right} cocycles, so $\beta'(g) = \beta(g^{-1})$ was also considered in \cite{JMP1}. A more effective approach in the same line was found in \cite{JMP2}, where only the left orthogonal cocycle $\beta$ and the lift $m = \widetilde{m} \circ \beta$ were required.

\begin{thm} \emph{(HM criterion in group algebras)}\textbf{.} \label{ThmHMCocycle} Let $\G$ be a locally compact unimodular group and $m \colon \mathrm{G} \to \C$. Let $\beta: \G \to \R^n$ be a cocycle map associated to an orthogonal action $\alpha \colon \G \curvearrowright \R^n$. Then, the following inequality holds for $1 < p < \infty$ and any lifting multiplier $m = \widetilde{m} \circ \beta$ $$\big\| T_m \colon L_p(\cL \G) \to L_p(\cL \G) \big\|_{\mathrm{cb}} \le C_p \sum_{|\gamma| \le [\frac{n}{2}] + 1} \Big\| |\xi|^{|\gamma|} \partial_\xi^\gamma \widetilde{m}(\xi) \Big\|_\infty.$$ 
\end{thm}

As in Remark \ref{Rem-Riesz-Schur}, the key idea to avoid noncommutative CZ methods was to express HM multipliers as Littlewood-Paley averages of Riesz transforms associated to fractional laplacians, at the cost of a worse constant $C_p$. This reduces the problem to dimension-free estimates for Riesz transforms in group algebras, in line with a vast commutative literature \cite{B1,DRdF,GV,LP3,Mey,PRiesz,St3}. The maps 
$$R_{\beta,u}(f) = \int_\G \frac{\langle \beta(g),u \rangle_\H}{\|\beta(g)\|_\H} \, \widehat{f}(g) \lambda_g \, d\mu(g) \quad \mbox{for each} \quad u \in \H$$ are noncommutative Riesz transforms on $\cL \G$ in terms of a cocycle map $\beta \colon \G \to \H$ associated to an orthogonal action $\alpha \colon \G \curvearrowright \H$. This includes infinite-dimensional Hilbert spaces $\H$, which are needed for Riesz transforms generated by fractional laplacians |unreachable with previous methods| and crucial in turn for Theorem \ref{ThmHMCocycle}. Dimension-free estimates are derived in \cite{JMP2} from probabilistic tools.  

\begin{rem} \label{Rem-Cocycles}
\emph{In Euclidean spaces, Theorem \ref{ThmHMCocycle} also recovers both de Leeuw's restriction and periodization theorems \cite{dL} and some exotic multipliers \cite[Section 5]{JMP1}. Links with noncommutative geometry appear in \cite{AK} and \cite[Appendix C]{JMP2}.} 
\end{rem}

The results in \cite{JMP2} recently inspired new results for Schur multipliers extending the original ones for Fourier multipliers. More precisely, by the last assertion in Theorem \ref{ThmFS2} i), the dimension-free estimates above can be transferred to their Schur analogues. Their nonToeplitz extensions for arbitrary Schur multipliers have been investigated by Arhancet and Krieger \cite[Theorem 3.3]{AK} and recently \cite{GPPR} in further detail. The approach in \cite{AK} closely follows \cite{JMP2}, while the later one gives a cleaner statement and prominently a much simpler proof. None of the usual analytic or probabilistic methods —Fourier transforms and CZ techniques or diffusion/Markov semigroups and Pisier’s reduction formula— are needed. On the contrary, the simpler argument in \cite{GPPR} is modeled on Grothendieck’s inequality \cite[Chapter 5]{PisSim} and the link below with Grothendieck's work is not accidental.

\begin{cor} \label{Cor-Gro}
Let $\Gamma$ be any index set, consider $\{u_j, u_j', w_j, w_j' \colon j \in \Gamma\}$ arbitrary vectors in a Hilbert space $\H$ and let $\Lambda: \H \to \H$ be a contraction. Then, the symbols $$M(j,k) = \Big\langle \frac{u_j + u_k'}{\|u_j + u_k'\|}, \Lambda \Big( \frac{w_j + w_k'}{\|w_j + w_k'\|}\Big)\Big\rangle$$ yield completely bounded Schur multipliers $S_M \colon S_p(\Gamma) \to S_p(\Gamma)$ for any $1 < p < \infty$.
\end{cor}

Taking $u_k' = w_j = 0$, we obtain symbols of the form $\langle \xi_j, \psi_k \rangle$ for some uniformly bounded families of vectors in $\H$. These characterize $S_\infty$-bounded Schur multipliers \cite{Gr,PisSim,PisBAMS} and Corollary \ref{Cor-Gro} gives a weaker form of the Grothendieck-Haagerup's criterion for Schatten $p$-classes. The classical boundedness of triangular truncations also follows taking $\Gamma = \R$, $(u_j, u_k', w_j, w_k') = (j,-k,1,0)$ and $\Lambda = \mathrm{id}$.

\begin{rem}
\emph{A H\"ormander-Mikhlin multiplier theory valid for free groups was recently established by Tao Mei, Éric Ricard and Quanhua Xu in \cite{MRX} with methods relying on Mei-Ricard's remarkable work \cite{MR} and very different from the ones above.}
\end{rem}

\subsubsection{H\"ormander-Mikhlin criteria via Lie derivatives} 

The cocycle approach above was inspired by the lack of differential structures on general topological groups. Lie groups though carry their own differential structure. Can we find more intrinsic H\"ormander-Mikhlin conditions for Lie groups? This is specially interesting for Lie groups lacking finite-dimensional orthogonal cocycles, as noncompact simple Lie groups. In fact, the results below were originally motivated by Lafforgue/de la Salle's theorem \cite{LdlS}, which we shall review in Section \ref{sec:rigidity}.   

Given a unimodular $n$-dimensional Lie group $\G$, consider the left-invariant vector fields in $\G$ generated by any given orthonormal basis $\mathrm{X}_1, \mathrm{X}_2, \ldots,  \mathrm{X}_n$ of its Lie algebra $\mathfrak{g}$. The left-invariant Lie derivatives $$\partial_{\mathrm{X}_j} m (g) \, = \, \frac{d}{ds}\Big|_{s=0} m \big( g \exp(s \mathrm{X}_j) \big)$$ do not commute for $j \neq k$. This justifies to define the set of multi-indices $\gamma$ as ordered tuples $\gamma = (j_1, j_2, \ldots, j_k)$ with $1 \le j_i \le \dim \G$ and $|\gamma| = k \ge 0$, which correspond to the Lie differential operators $$d_g^\gamma m (g) \, = \, \partial_{\mathrm{X}_{j_1}} \partial_{\mathrm{X}_{j_2}} \cdots \, \partial_{\mathrm{X}_{j_{|\gamma|}}} m (g).$$
We start with a local form of the H\"ormander-Mikhlin theorem for Lie groups \cite{CGPT2}.

\begin{thm} \label{Thm-LocalHM} \emph{(Local HM criterion)}\textbf{.} Consider a finite-dimensional unimodular Lie group $\G$ with Riemannian metric $\rho$ and set $L_\mathrm{R}(g) =\rho(g,e)$. Let $m \colon \G \to \C$ be supported by a relatively compact neighborhood of the identity $\Omega$. Then, the following inequality holds for $1 < p < \infty$ $$\big\| T_m \colon L_p(\cL \G) \to L_p(\cL \G) \big\|_{\mathrm{cb}} \, \le \, C_p(\Omega) \sum_{|\gamma| \le [\frac{\dim \G}{2}]+1} \big\| L_\mathrm{R}(g)^{|\gamma|} d_g^\gamma m(g) \big\|_\infty.$$
\end{thm} 

This intrinsic H\"ormander-Mikhlin condition for Lie groups has sharp regularity orders and is necessarily local, since its validity for arbitrary symbols would get in conflict with \cite{LdlS}. Although a highly technical argument was earlier presented in \cite{PRS} for special linear groups, Theorem \ref{Thm-LocalHM} follows from Theorem \ref{ThmHMS}:

i) \textbf{Local Fourier-Schur transference.} We first use Theorem \ref{ThmFS2} ii) to rewrite the problem in terms of Herz-Schur multipliers. In fact, by relative compactness of $\Omega$ and smooth partitioning, we may assume that $\Omega$ is as small as we need. On the other hand, given $\Lambda$ any neighborhood of the identity in $\G$ and by a refined version of local transference \cite{PRS,PST}, we may also assume that $M_1(g,h) = m(gh^{-1})$ is just defined over $\Lambda \times \Lambda$ as long as $\Omega$ is small enough. 

ii) \textbf{Local lifting into the Lie algebra.} Let $\mathfrak{g}$ be the Lie algebra of $\G$. Since the exponential map $\exp \colon \mathfrak{g} \to \G$ is a local diffeomorphism at the identity, we may easily construct a new symbol $M_2 \colon \mathfrak{g} \times \mathfrak{g} \to \C$ so that $M_2(x,y) = m(\exp(x) \exp(y)^{-1})$ for $x,y \in \exp^{-1}(\Lambda)$ with $\Lambda$ small enough and such that $M_2$ is smooth away from the diagonal and compactly supported. By \cite[Lemma 2.1]{PST} both $S_{M_1}$ and $S_{M_2}$ have the same cb-norm on Schatten $p$-classes, which opens a door to Theorem \ref{ThmHMS}.

This reduces a Fourier multiplier problem to a nonToeplitz Schur multiplier one. All what is left is to relate Lie and Euclidean metrics/derivatives, which are locally equivalent near the identity. This whole argument is an illustration of the strength of nonToeplitz harmonic analysis and the great flexibility of Schur multipliers to preserve their $S_p$-mapping bounds under manipulations. 
    
\null

\vskip-20pt

\null

\begin{center}
\includegraphics[scale=0.8]{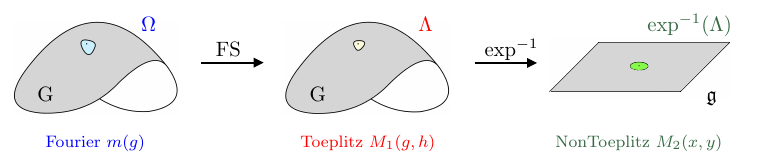}
\end{center}

\vskip-20pt

\null

\begin{center}
{\scriptsize \textbf{Figure 4.} Local lifting of Fourier multipliers to nonToeplitz Schur multipliers \\ [-3pt] This has further applications in harmonic analysis over Lie groups, see Section \ref{Subsect-HilbertLie}} 
\end{center}

\vskip5pt

Motivated once more by \cite{LdlS}, a natural goal is to eliminate locality for simple Lie groups. The natural length $L_\G \colon \G \to \R_+$ for this class of groups is locally Euclidean around the identity and its asymptotic behavior is dictated by the adjoint representation $L_\G(g) \approx \|\mathrm{Ad}_g\|^{\tau_\G}$ as $g \to \infty$ for $\tau_\G = \mathrm{d}_\G/[\frac{\dim \G + 1}{2}]$, with $\mathrm{d}_\G$ from \cite{Mauc}. Maucourant's constant $\mathrm{d}_\G$ gives the volume growth rate of $\mathrm{Ad}$-balls up to a logarithmic factor. It turns out that $(\tau_\G, \mathrm{d}_\G) = (1/2, (n^2-1)/4)$ for $\G = S \hskip-1pt L_n(\R)$ and $\tau_\G \le 1$ for any simple Lie group.

\begin{thm} \label{Thm-GlobalHM} \emph{(Global HM for simple Lie groups)}\textbf{.} 
Let $\G$ be a finite-dimensional simple Lie group with $\mathrm{d}_\G \ge 2 \mbox{$[\frac12 (\dim \G + 1)]$}/\dim \G$. Then, the following inequality holds for any Fourier symbol $m \colon \G \to \C$ and $1 < p < \infty$ $$\big\| T_m \colon L_p(\cL \G) \to L_p(\cL \G) \big\|_{\mathrm{cb}} \, \le \, C_p \sum_{|\gamma| \le [\frac{\dim \G}{2}]+1} \big\| L_\G(g)^{|\gamma|} d_g^\gamma m(g) \big\|_\infty.$$
\end{thm} 

This result from \cite{CGPT2} significantly improves \cite[Theorem A]{PRS}. The assumption $\mathrm{d}_\G \ge 2 \mbox{$[\frac12 (\dim \G + 1)]$}/\dim \G$ ensures that the H\"ormander-Mikhlin condition above implies asymptotically $$|d_g^\gamma m(g)| \lesssim \|\mathrm{Ad}_g\|^{- \mathrm{d}_\G} \quad \mbox{for} \quad \mbox{$|\gamma| \le [\frac{n}{2}]+1$.}$$ 

\begin{rem}
\emph{This decay for $\gamma=0$ is in line with Lafforgue/de la Salle's theorem \cite{dLdlS,LdlS}. In fact, refined necessary conditions \cite{PRS} for radial multipliers $m(g) = \varphi(|g|)$ in terms of the Hilbert-Schmidt norm $| \ |$ also show a decay at infinity for a number of derivatives of $\varphi$ as long as $T_m$ is $L_p$-bounded and $p$ is large enough. This necessity increases with the rank and there exists radial multipliers satisfying Theorem \ref{Thm-GlobalHM} in a given rank $n$ and failing the necessary conditions for ranks $m >> n$.}
\end{rem}

\begin{rem}
\emph{The optimal regularity order in Theorem \ref{Thm-GlobalHM} leads to the critical decay order $\mathrm{d}_\G$ for $m$. As noted in \cite[Remark 3.2]{CGPT2}, there are some indications that there is no more room for improvement in the metric $L_\G$. Also, condition $\mathrm{d}_\G \ge 2 \mbox{$[\frac12 (\dim \G + 1)]$}/\dim \G$ holds for large classes of simple Lie groups but fails for $S \hskip-1pt L_2(\R)$, which is weakly amenable. Thus, it is expectable but open to find HM conditions with arbitrarily slow decay for $S \hskip-1pt L_2(\R)$ and other rank one simple Lie groups, beyond the scope of Theorem \ref{Thm-GlobalHM}. In this direction, Martijn Caspers has recently obtained in \cite{C1} lower asymptotic decay rates for a class of K-biinvariant smooth symbols, which fail to admit a singular behavior around the identity. This interesting result resembles Calderón–Torchinsky interpolation theorem \cite{CT}.}
\end{rem}

\subsection{Characterizing Hilbert transforms on Lie groups} \label{Subsect-HilbertLie} 

Idempotent Fourier multipliers are those whose symbols are the characteristic function of certain domain $\Omega$. That is, Fourier truncations over the frequencies lying on $\Omega$. Historically, these multipliers have been considered in the problem of $L_p$-convergence for Fourier series and integrals. In Euclidean spaces, iterating half-space multipliers along several directions (Hilbert transforms) shows that convex polyhedra are valid examples of Fourier $L_p$-idempotents for $1 < p \neq 2 < \infty$, while Fefferman's ball multiplier theorem \cite{Fe} confirms that the boundary $\partial \Omega$ must indeed be flat a.e. This led additionally to new connections with highly singular operators coming from Kakeya sets or Bochner-Riesz means, central in Euclidean harmonic analysis since then. 

The same problem has been studied for other groups. Which domains $\Omega$ of a given group $\G$ define a Fourier idempotent in $L_p(\cL \G)$? In this case, both the shape of $\Omega$ and the geometry of $\G$ play a key role. Mei-Ricard's fundamental work on free groups \cite{MR} along with \cite{PRo} were the first contributions. The local geometry of Fourier $L_p$-idempotents on Lie groups is particularly interesting, since the curved geometry of general Lie groups makes unclear what should be the notion of \lq\lq boundary flatness\rq\rq${}$ on them. Once more, Schur multipliers come into play to solve this problem. Indeed, Theorem \ref{ThmIdemp} has led to a characterization of the local geometry of Fourier $L_p$-idempotents on Lie groups, 
applying the process sketched in Figure 4. Before stating this result, let us consider three fundamental examples of a group $\G$ with a smooth domain $\Omega$:
\begin{itemize}
\item[$\bullet$] The real line $\G_1 = \R$ with $\Omega_1 = (0,\infty)$. \vskip2pt

\renewcommand{\thefootnote}{$\star$}

\item[$\bullet$] The affine group $\G_2 = \Aff_+(\R)$\footnote{Affine increasing bijections $x\mapsto ax+b$ for $a \in \R_+^*$ and $b \in \R$, isomorphic to $\R \rtimes \R_+^*$.} and $\Omega_2 = \{ax+b: b > 0\}$.

\renewcommand{\thefootnote}{$\dag$}

\item[$\bullet$] The universal covering group $\G_3 = \widetilde{P \hskip-1pt S \hskip-1pt L}_2(\R)$\footnote{The action $\alpha: \widetilde{P \hskip-1pt S \hskip-1pt L}_2(\R) \curvearrowright \R$ is obtained by lifting the standard action of $P \hskip-1pt S \hskip-1pt L_2(\R)$ on the projective line to the universal covers. If $p \colon \R \to P^1(\R)$ denotes the universal cover, then the universal cover of $\SL_2(\R)$ is identified with the group of homeomorphisms $g: \R \to \R$ for which there is $A \in P \hskip-1pt S \hskip-1pt L_2(\R)$ such that $p\circ g= A \cdot p$.} with $\Omega_3 = \{g: \alpha_g(0) > 0\}$.
\end{itemize}
These Hilbert transforms are (globally) Fourier cb-$L_p$-idempotents, $1 < p < \infty$ \cite{GPX}. 
Also, $m \colon \G \to \C$ defines a locally bounded Fourier $L_p$-multiplier at $g_0 \in \G$ when there is a function $\varphi \equiv 1$ on a neighbourhood of $g_0$ and $T_{\varphi m} \colon L_p(\cL \G) \to L_p(\cL \G)$.

\begin{thm} \label{DDD} \emph{(Geometry of Fourier idempotents on Lie groups)}\textbf{.} \label{Thm-IdempLie} Let $\G$ be a simply connected Lie group. Consider  a $\mathcal{C}^1$-domain $\Omega$ in $\G$ and let $g_0 \in \partial \Omega$. Then the following are equivalent for $1 < p \neq 2 < \infty$\emph{:} 
\begin{itemize}
\item[i)] $\chi_\Omega$ defines locally at $g_0$ a completely bounded Fourier $L_p$-multiplier.

\item[ii)] $\partial \Omega = g_0 \exp(\mathfrak{h})$ locally near $g_0$ for some codimension $1$ Lie subalgebra $\mathfrak{h}$.

\item[iii)] There exists a smooth surjective homomorphism $\zeta \colon \G \to \G_j$ for some index $j \hskip-3pt = \hskip-3pt 1,2,3$ for which the identity $\Omega \hskip-1.5pt = \hskip-1.5pt g_0 \zeta^{-1}(\Omega_j)$ holds in a neighborhood of $g_0$
\end{itemize}
\end{thm}


Theorem \ref{Thm-IdempLie} was established in \cite{PST} and condition ii) makes clear what \lq\lq Fourier boundary flatness\rq\rq${}$ means for Lie groups. A reformulation in terms of Lie group actions by diffeomorphisms on $\R$ |classified by Lie himself \cite{zbMATH02681669}| then leads to condition iii). Unexpectedly, all Fourier idempotents arise from the \emph{classical}, the \emph{affine} and the \emph{projective} Hilbert transforms above. In Euclidean spaces, Theorem \ref{Thm-IdempLie} recovers that Fourier $L_p$-bounded idempotents locally correspond to half-space multipliers, which are directional amplifications of the Hilbert transform on $\R$. By analogy, Fourier cb-$L_p$-idempotents on arbitrary Lie groups arise as \emph{directional amplifications} of one of these three fundamental models.

\begin{rem}
\emph{Theorem \ref{Thm-IdempLie} gives a clear description of all Hilbert transforms on nilpotent Lie groups. It also shows that simple Lie groups lack to admit completely $L_p$-bounded Fourier idempotents for $p \neq 2$, except those locally isomorphic to $S \hskip-1pt L_2(\R)$, which carry a unique local Fourier idempotent up to left/right translations.}
\end{rem}

\begin{rem}
\emph{A natural problem for Fourier idempotents over discrete groups is to study the $L_p$-convergence of a sequence of compact Fourier truncations. In this framework, groups like $S \hskip-1pt L_2(\Z)$ or free groups are of special interest. Being weakly amenable they admit smooth Fourier approximations, but further insight on more singular Fourier approximations could uncover structural properties of their group algebras. A renowned challenge is the \emph{free ball multiplier problem} on the (failure) of uniform $L_p$-bounds for Fourier truncations over the Cayley graph balls. This is well-known \cite{BF,Ku} for $|\frac12 - \frac1p| \ge \frac16$ but remains open for values of $p$ closer to $2$.}
\end{rem}

\section{\bf \large Operator rigidity phenomena in higher-rank lattices}\label{sec:rigidity}

The structure and classification of group von Neumann algebras and their crossed products with probability spaces |group measure spaces| is a fundamental topic in operator algebras. How much information about a discrete group $\Gamma$ is retained in $\mathcal{L} \Gamma$ or its measure spaces $L_\infty(\Omega) \rtimes \Gamma$? Connes' celebrated classification of injective factors \cite{Co1} implies that discrete amenable groups and their actions give rise to undistinguishable factors. In sharp contrast, group algebras of nonamenable groups may even remember the complete structure of the group. Therefore, group von Neumann algebras range from strikingly nonrigid to extremely rigid ones and the classification of nonamenable II$_1$ factors remains largely intractable. 

In 1980, Connes conjectured \cite{CRC} that $\mathcal{L}\Gamma \simeq \mathcal{L}\Lambda \Rightarrow \Gamma \simeq \Lambda$ for any pair $(\Gamma,\Lambda)$ of property (T) i.c.c. discrete groups. This bold assertion  somehow claims that these groups are \lq\lq nonamenable enough\rq\rq{} to be pairwise distinguishable from their group algebras. A key instance of Connes' conjecture  refers to $\Gamma, \Lambda \in \{P\hskip-1pt S \hskip-1pt L_n(\Z) \hskip-2pt: n \ge 3\}$ or other higher rank lattices. We refer to \cite{Hou,dlS} for recent discussions of this conjecture. Rigidity theory for von Neumann algebras (operator rigidity in what follows) is a challenging subject. Popa's deformation/rigidity theory and impressive findings thereafter \cite{BH,HSV,Io,Oz1,OP1,OP2,Pe, Po1, Po2, Po3, PV1, PV2} 
illustrate a great progress for group measure spaces, but the superrigidity problem for group factors remains rather incomplete|Ioana/Popa/Vaes found in 2013 the first examples \cite{IPV} and the first property (T) ones were recently found in \cite{CIOS}. 

\subsection{Lafforgue/de Laat/de la Salle's rigidity theorem}

In what follows, an approximate identity (AI) is a sequence $\Phi = \{\phi_j \colon j \ge 1\} \subset \mathcal{C}_\mathrm{c}(\G)$ converging to $1$ uniformly over compacta. Leptin's well-known characterization claims that $\G$ is amenable if and only if there exists an approximate identity $\Phi$ made of positive definite functions $\phi_j$. The operator cb-norm of the Herz-Schur multipliers $S_{\phi_j}$ is then $1$. $\G$ is called weakly amenable or is said to have the Cowling/Haagerup approximation property (AP) when there exists an AI $\Psi = \{ \psi_j \colon j \ge 1\}$ of (not necessarily positive definite) functions such that the Schur multipliers $\{S_{\psi_j} \colon j \ge 1\}$ are uniformly bounded in the operator cb-norm. $\G$ is \emph{Schur weakly $p$-amenable} for some $p > 2$ when instead $$\mathrm{Sch}_p(\G) := \inf_{\Psi \, \mathrm{AI}} \, \sup_{j \ge 1} \big\| S_{\psi_j} \colon S_p(\G) \to S_p(\G) \big\|_{\mathrm{cb}} < \infty.$$ This was introduced in \cite{LdlS} as Schur AP (approximation property) of $S_p(\G)$, but we prefer to refer to it as a group property. It is clear that $\mathrm{Sch}_2(\G) = 1$ for every locally compact group. Moreover, by an elementary interpolation argument, Schur weak $p$-amenability becomes increasingly restrictive as $p \to \infty$. The case $p = \infty$ corresponds to weak amenability, and $\mathrm{Sch}_\infty(\G)$ is known as the Cowling/Haagerup constant of $\G$|a von Neumann algebra invariant retained by $\cL \G$. In sum, the Schur critical index $q_\G$ below which $\G$ is Schur weakly $p$-amenable \emph{measures the degree of nonamenability} of $\G$. Lafforgue/de la Salle's work, along with subsequent results with de Laat, implies that higher rank simple Lie groups and their lattices can be as nonamenable as we want |within this scale| as the rank increases. In the key case of $S \hskip-1pt L_n(\R)$ and $S \hskip-1pt L_n(\Z)$, we compile \cite{dLdlS,LdlS} as follows. 

\begin{thm} \label{Thm-LdlSRigidity} \emph{(Lafforgue-de Laat-de la Salle's rigidity)}\textbf{.} Both $S \hskip-1pt L_n(\R)$ and $S \hskip-1pt L_n(\Z)$ fail to be Schur weakly $p$-amenable for $p > 2+ \frac{2}{\alpha_n}$, with $\alpha_n = \lfloor \frac{n-1}{2} \rfloor \to \infty$ as $n \to \infty$.
\end{thm}

This strengthens a celebrated result by Uffe Haagerup, who disproved the weak amenability for higher rank simple Lie groups \cite{CDSW,H1}. Theorem \ref{Thm-LdlSRigidity} applies to any lattice and implies |by the last assertion in Theorem \ref{ThmFS2} i)| the failure of similar $L_p$-approximation properties for Fourier multipliers over $S \hskip-1pt L_n(\R)$ and $S \hskip-1pt L_n(\Z)$ for the same values of $p$, see below. This has a double impact: 

$\bullet$ In harmonic analytic terms, it gives the first class of groups over which no Fourier $L_p$-approximation works |regardless of how smooth the (compact) Fourier truncations are| for finite values of $p$. This highlights dramatic pathologies in harmonic analysis over simple Lie groups and lattices. 

$\bullet$ In operator algebraic terms, this unprecedented phenomenon opens a new door to attack Connes' rigidity conjecture. Indeed, it would suffice to construct Fourier $L_p$-approximations over $P \hskip-1pt S \hskip-1pt L_n(\Z)$ when $p$ is close enough to $2$, since the critical index for which this happens is a von Neumann algebra invariant. 

\noindent This shows why Lafforgue/de Laat/de la Salle's theorem stands as a rigidity result.  

\begin{rem} \label{Rem-LieLattice} 
\emph{$\mbox{Sch}_p(\G) = \mbox{Sch}_p(\Gamma)$ for any lattice $\Gamma$ in $\G$. This follows by good restriction/extension properties of Schur multipliers \cite[Theorem 2.5]{LdlS}. Restriction of Fourier multipliers in group algebras was investigated in \cite{CPPR}. Unfortunately, no general restriction theorem seem to be within reach with current techniques.}
\end{rem}

\subsection{On Fourier and Schur $L_p$-approximations} \label{Subsect-FSAP}

Given $p>2$, we say that $\G$ is \emph{Fourier weakly $p$-amenable} when there exists an approximate identity $\Psi$ whose Fourier multipliers $T_{\psi_j}$ remain uniformly cb-bounded in $L_p(\cL \G)$. That is, we just replace Schur multipliers above by their Fourier peers. Theorems \ref{ThmFS2} and \ref{Thm-LdlSRigidity} imply that $S \hskip-1pt L_n(\Z)$ fails Fourier weak $p$-amenability for $p > 2 + 2/\alpha_n$. By \cite{JR}, this can be reformulated by saying that $L_p(\cL \G)$ fails the \emph{completely bounded approximation property} (CBAP -- see e.g. \cite{CH,H,H1}) for those values of $p$. The CBAP is a von Neumann algebra invariant and this has potential consequences in Connes' rigidity problem for $P \hskip-1pt S \hskip-1pt L_n(\Z)$. Indeed, in case one could justify that $S \hskip-1pt L_{2n-1}(\Z)$ is Fourier weakly $p$-amenable for some $p > 2$, then we would have $$\cL P \hskip-1pt S \hskip-1pt L_{2n-1}(\Z) \neq \cL P \hskip-1pt S \hskip-1pt L_{2m-1}(\Z) \quad \mbox{for} \quad m > \frac{p}{p-2}.$$ Upper/lower bounds for the Fourier critical weak $p$-amenability index of $S\hskip-1pt L_n(\Z)$ would refine Theorem \ref{Thm-LdlSRigidity} and the above implication on Connes' rigidity problem.

\begin{pro} Is the Fourier critical index $p_n$ greater than $2$ for every $n \ge 3$\emph{?} This would solve Connes' rigidity $\mathcal{L} P \hskip-1pt S \hskip-1pt L_n(\Z) \neq \mathcal{L} P \hskip-1pt S \hskip-1pt L_m(\Z)$ at least for infinitely many pairs $m,n \ge 3$. More generally, is there any \textbf{change of regime} for $L_p(\mathcal{L} P \hskip-1pt S \hskip-1pt L_n(\Z))$ at certain critical index $\widetilde{p}_n > 2$ which is retained by the group von Neumann algebra of $P \hskip-1pt S \hskip-1pt L_n(\Z)$ and simultaneously satisfies that $\widetilde{p}_n \to 2$ as $n \to \infty$\emph{?}
\end{pro}

\begin{rem} 
\emph{A related intriguing question going back to Lafforgue/de la Salle's work \cite{LdlS} is whether or not Fourier and Schur weak $p$-amenability are equivalent properties. They are for amenable groups and Fourier weak $p$-amenability implies its Schur analog for nonamenable ones, as it follows from Theorem \ref{ThmFS2}. Thus, in the terminology from \cite{LdlS}, we still ignore whether the Schur approximation property for $S_p(\G)$ implies the CBAP for $L_p(\cL \G)$. No expert would bet a lot in a global form of Fourier-Schur transference for nonamenable groups |no counterexamples are known either| but the question above is definitely wider.}
\end{rem} 

Is Schur weak amenability a von Neumann algebra invariant for discrete groups? This is an even wider formulation of the above question. Besides some strategies to address this problem, the interest in Schur weak amenability goes beyond. Showing nontrivial Schur weak amenability indices or even finding weaker changes of regime for $S \hskip-1pt L_n(\R)$ or $S \hskip-1pt L_n(\Z)$ |to be discussed in Section \ref{Sub-Radial-Polinomial}| would not only be of independent interest, but also a model to follow for Fourier approximations, in case no general relation between both can be established. In particular, we now present some results and comments on Schur multipliers pointing in this direction.  

According to Remark \ref{Rem-LieLattice} we may work in $S \hskip-1pt L_n(\R)$ instead of $S \hskip-1pt L_n(\Z)$, where the underlying topology and the additional structure from the Iwasawa decomposition could help in relating the problem with somehow modern techniques from Euclidean harmonic analysis. Back in 2015, we conjectured a link between Schur critical indices with Bochner-Riesz critical exponents.
\begin{conj} \label{Conj-SCI}
\emph{The Schur critical indices of simple Lie groups/lattices satisfy \\ $q_\G := \sup \big\{ p \ge 2 \colon \G \mbox{ is Schur weakly $p$-amenable} \big\}$ equals $2 \, \mathrm{rank}(\G) / \mathrm{rank}(\G) -1$}.
\end{conj}
This holds for $S \hskip-1pt L_2(\R)$, matches Theorem \ref{Thm-LdlSRigidity} for $S \hskip-1pt L_3(\R)$ and refines it for $n \ge 3$. Its Fourier analog for simple lattices $\Gamma \subset \G$ implies that $\mathrm{rank}(\G)$ is retained in $\cL \Gamma$|fully solving Connes' rigidity for $\{P \hskip-1pt S \hskip-1pt L_n(\Z) \colon n \ge 3\}$. Is also $q_\G$ a von Neumann algebra invariant? If not, can we prove a Fourier analog? In Section \ref{Sub-Radial-Polinomial} below, we consider a less abrupt \emph{change of regime} |at the same critical index determined by the rank| which is formally simpler, but still potentially useful to face Connes' rigidity if Conjecture \ref{Conj-SCI} became false or inaccessible. This conjecture was originally based on the behavior of Fourier multipliers $T_\phi$ for smooth compactly supported $\phi \colon S \hskip-1pt L_n(\R) \to \R_+$. As shown in \cite{PRo}, these multipliers are mirrored in their lifts $T_{\Phi} \rtimes \mathrm{id}$ over $\R^n \rtimes S \hskip-1pt L_n(\R)$, but the action $S \hskip-1pt L_n(\R) \curvearrowright \R^n$ produces a severe loss of regularity in $T_\phi$ when $\mathrm{supp} \, \phi$ is large, since large elements $g \in S \hskip-1pt L_n(\R)$ combine rotations and volume-preserving dilations with arbitrarily high eccentricity. Again inspired by \cite{PRo}, $L_p$-boundedness for $T_{\Phi} \rtimes \mathrm{id}$ imposes a \emph{uniform control} of these actions. Bochner-Riesz multipliers are the closest model of Euclidean $L_p$-multipliers for $|\frac12 - \frac1p|$ small enough in terms of the dimension. Not in vain Conjecture \ref{Conj-SCI} gives the Bochner-Riesz exponents in dimension $\mathrm{rank} (\G)$.  

\renewcommand{\theequation}{RdF}

Given $\phi_\delta(\xi) = (1 - |\xi|^2)_+^\delta$ for $\delta >0$, the Bochner-Riesz conjecture claims that the Fourier multiplier $T_{\phi_\delta}$ is $L_p$-bounded on $\R^n$ iff $|\frac12 - \frac1p| < \frac{1+\delta}{2n}$. It has only been confirmed in dimension 2 and stands as one of the hardest problems in harmonic analysis. It is also discouraging that the $\mathcal{B}(L_p)$-norm of $T_{\phi_\delta}$ is arbitrarily large as $\delta \to 0$ for any $p\neq 2$. The above analogy with rigidity thus suggests that no uniform $L_p$-bounds are possible for a family of multipliers in $S \hskip-1pt L_n(\R)$ with arbitrarily large supports. Interestingly, the work of Córdoba and Rubio de Francia \cite{Cor,RdF} leads to uniform bounds  
\begin{equation} \label{EqRdF}
\sup_{\delta>0} \big\| T_{\phi_\delta} \colon \Lambda_{p2}(\R^n) \to \Lambda_{p2}(\R^n) \big\| < \infty \ \Leftrightarrow \ \frac{2n}{n+1+2\delta} < p < \frac{2n}{n-1-2\delta},
\end{equation}
where $\Lambda_{p2}(\R^n)$ is the mixed $L_p(L_2)$-norm space in polar coordinates, with $L_p$-norm for the radial variable and $L_2$-norm for the angular ones. In view of Fefferman's theorem for the ball \cite{Fe}, using an $L_2$-norm in the variables generating curvature explains this better behavior. Some time ago, we introduced a mixed-norm space $\Lambda_{p2}$ in the matrix-algebra of simple Lie groups, with \lq polar coordinates\rq{} coming from the Iwasawa decomposition $\G = \mathrm{KS} = \mathrm{SK}$ 
$$\|A\|_{\Lambda_{p2}(\G)} = \Big\| \Big( \int_\mathrm{K} (A^*A)_{gk,hk} \, dk \Big)^\frac12 \Big\|_{S_p(\G)}.$$ 
A priori, applications towards Connes' rigidity would require $\Lambda_{p2}$-extensions of Theorem \ref{Thm-LdlSRigidity}. In particular, it was certainly shocking to find that, for simple Lie groups, the \textbf{Schur approximation property is equivalent to its mixed-norm analogue}. Indeed, if $\phi \colon \G \to \C$ is K-biinvariant \renewcommand{\theequation}{MN$_p$}
\begin{equation} \label{Eq-Mixed-Norm}
\big\| S_\phi \colon \Lambda_{p2}(\G) \to \Lambda_{p2}(\G) \big\|_{\mathrm{cb}} = \big\| S_\phi \colon S_p(\G) \to S_p(\G) \big\|_{\mathrm{cb}}.
\end{equation}
This unexpected property fails in Euclidean spaces|$L_p$-bounds do not compare to mixed-norms since polar decomposition $\R^2 = \mathbf{S}^1 \times \R_+$ lacks the group structure in Iwasawa's decomposition. As a consequence of this and a matrix-valued form of a well-known duality argument in harmonic analysis, we also found: 

\begin{thm} \label{ThmKakeya} \emph{(A sufficient condition for the Schur AP)}\textbf{.} Consider the normalized gaussians $\gamma_q(x) = \gamma(x)/\|\gamma\|_q$ for $\gamma(x) = \exp(-|x|^2)$. Then, $S \hskip-1pt L_n(\R)$ is Schur weakly $p$-amenable provided the inequality below holds for $q = \frac{p}{p-2}$, every \emph{K}-Toeplitz matrix $A$ $(A_{gk,hk} = A_{g,h}$ for all $k \in \mathrm{K})$ and any unit vector $u \in \R^n$ $$\Big( \int_{\R^n} \Big\| \sup_{\begin{subarray}{c} |w|=1 \\ L > 0 \end{subarray}} - \hskip-10.8pt \int_{-L}^L \gamma_q(x + sw) \Big( e^{2\pi i s \langle (g^{-1} - h^{-1})u,w \rangle} A_{gh} \Big) ds \Big\|_{S_q(\G)}^q dx \Big)^\frac1q \le_{\mathrm{cb}} \|A\|_{S_q(\G)}.$$ Its validity for any $q > n$ readily implies Schur weak $p$-amenability for $2 \le p < \frac{2n}{n-1}$.
\end{thm}

This sufficient condition is a noncommutative maximal inequality generalizing the Kakeya universal maximal operator to a class of matrix-valued functions. This maximal is unbounded in Euclidean $L_q$-spaces, but turns out to be $L_q(\R^n)$-bounded for $q > n$ when acting on radial functions. It is precisely \eqref{Eq-Mixed-Norm} what allows us to treat the matrix $A$ as K-Toeplitz, a matrix-valued form of radiality. This somehow supports the conjectured connection between Connes' rigidity and Bochner-Riesz multipliers / Kakeya maximal functions, it will be the subject of a forthcoming work. These harmonic analytic methods are yet to be explored in this context. 

\subsection{Weaker changes of regime -- Polynomial decay} \label{Sub-Radial-Polinomial}

Theorem \ref{Thm-LdlSRigidity} admits a quantitative form in terms of fast decay of K-biinvariant Schur multipliers. In the case of $S \hskip-1pt L_3(\R)$, if $d_s = \mathrm{diag}(e^s,1,e^{-s})$ and $m \in \mathcal{C}_0(S \hskip-1pt L_3(\R))$ is a $SO_3$-biinvariant symbol, the following holds for $p>4$  
\begin{equation} \tag{Exp$_p$} \label{EqExpDecay}
|m(d_s)| \le C_p e^{-c (1 - \frac{4}{p}) |s|} \big\| S_m \colon S_p(S \hskip-1pt L_3(\R)) \to S_p(S \hskip-1pt L_3(\R)) \big\|_{\mathrm{cb}},
\end{equation}
disproving the Schur AP. Given $M_\phi \colon \mathbf{S}^2 \times \mathbf{S}^2 \ni (x,y) \mapsto \phi(\langle x,y \rangle) \in \C$, the key point to prove the exponential decay \eqref{EqExpDecay} is to show that $\phi$ is $(\frac{1}{2} - \frac2p)$-H\"older when the Schur multiplier $S_{M_\phi}$ id cb-$S_p$-bounded. Similar considerations apply in rank $n$, with $n$-dimensional spheres and $p > \frac{2n}{n-1}$. This form of harmonic analysis on the sphere mirrors Euclidean phenomena|\emph{$\alpha_p$-H\"older regularity is necessary to define radial $L_p(\R^n)$-multipliers with $p$ above the BR exponent  $q_n = \frac{2n}{n-1}$}. An important change of regime is predicted by the Bochner-Riesz conjecture at its critical index $q_n$ |which provides radial Fourier $L_{q_n}$-multipliers in $\R^n$ failing to be $\alpha$-H\"older for arbitrarily small $\alpha$| and related changes of regime could take place in the Fourier and Schur multiplier theories over simple Lie groups.

This strongly motivates nonToeplitz forms of Bochner-Riesz multipliers \cite{PT} in 2D. Indeed, by product diffeomorphism stability \cite[Lemma 2.1]{PST}, this would imply Schatten $p$-class bounds for \emph{spherical Bochner-Riesz means} $\Phi_\delta(x,y) = \langle x,y \rangle_+^\delta$ over $\mathbf{S}^2 \times \mathbf{S}^2$ with $2 \le p \le 4$ and any $\delta>0$. Next, given any $\alpha>0$, we notice that $\Phi_\delta(x,y) = \phi(\langle x,y \rangle)$ and $\phi(\xi) = \xi_+^\delta$ is not $\alpha$-H\"older for $\delta$ small enough.  According to \cite{dLdlS}, this suggests that the decay in \eqref{EqExpDecay} might be  subexponential below the critical index. In fact, similar techniques yield $L_p$-bounded radial 2D multipliers with logarithmic modulus of continuity|like $\log(1 - |\xi|^2)^{-\beta_p}$ for large enough $\beta_p > 0$ and $2 \le p \le 4$, indicating that also polynomial rates of decay could be expected in $S \hskip-1pt L_3(\R)$. More generally, even if Conjecture \ref{Conj-SCI} or weaker forms of it were false or unverifiable, it seems that certain \lq\lq exponential $\to$ polynomial\rq\rq${}$ change of regime could hold for $p$ below $2 \, \mathrm{rank}(\G) / \mathrm{rank}(\G) -1$ in the context of \eqref{EqExpDecay}. 

In fact, one could go further and formulate the following weakening of Conjecture \ref{Conj-SCI}. In $S \hskip-1pt L_n(\R)$, set $\hskip-3pt \left\bracevert \hskip-2pt g \hskip-2pt \right\bracevert \hskip-3pt := \max \{\log \|g\|, \log \|g^{-1}\|\}$. Then, the question is whether there exists an approximate identity $\{\phi_j \colon j \ge 1\} \subset \mathcal{C}_\mathrm{c}(S \hskip-1pt L_n(\R))$ factorizing as $\phi_j (g) = m_j(g) (1 + \hskip-3pt \left\bracevert \hskip-2pt g \hskip-2pt \right\bracevert \hskip-3pt)^N$ for some $N \ge 0$ and satisfying $$\displaystyle \sup_{j \ge 1} \big\| S_{m_j} \colon S_p(S \hskip-1pt L_n(\R)) \to S_p(S \hskip-1pt L_n(\R)) \big\|_{\mathrm{cb}} < \infty \quad \mbox{for} \quad 2 \le p \le \frac{2(n-1)}{n-2}.$$ Similarly, given any other simple Lie group $\G$, one could set $\hskip-3pt \left\bracevert \hskip-2pt g \hskip-2pt \right\bracevert \hskip-3pt \approx \log (1+L_\G(g))$ for the length $L_\G \colon \G \to \R_+$ introduced before Theorem \ref{Thm-GlobalHM} and wonder about the same inequality for $2 \le p \le 2 \, \mathrm{rank}(\G) / \mathrm{rank}(\G) -1$. By Schur restriction, the above property implies the same assertion for any lattice $\Gamma \subset \G$.
This is a stronger change of regime and recovers Conjecture \ref{Conj-SCI} when $N=0$. It can also be regarded as a Schur AP over certain Sobolev spaces admitting $N$ derivatives in $L_p$. The ideas in Section \ref{Subsect-FSAP} could also be useful in this direction, while the allowed polynomial decay for $m_j$ should make it more accessible. In harmonic analytic terms, this result would be very interesting. In operator algebraic terms, it is not yet clear how to relate it with a vNa invariant of $\cL \Gamma$ |due to the length $\hskip-3pt \left\bracevert \hskip-2pt {\color{white}g} \hskip-2pt \right\bracevert \hskip-3pt$, not even its Fourier analog| but it seems conceivable that such an abrupt change of regime would be retained in the group algebra of higher rank simple lattices. 

\vskip5pt

\noindent \textbf{Acknowledgement.} \\ I would like to thank all my collaborators for so many exciting ideas, discussions and discoveries together in this beautiful subject... and also for the ones to come!! I also thank A. González-Pérez, M. de la Salle and E. Tablate for helpful comments.

\bibliographystyle{amsplain}

\begin{thebibliography}{100}

\bibitem {AK} C. Arhancet and C. Kriegler, Riesz transforms, Hodge-Dirac operators and functional calculus for multipliers. Lecture Notes in Math. \textbf{2304}. Springer, 2022.

\bibitem {B1} D. Bakry, Transformation de Riesz pour les semi-groupes sym\'etriques. S\'eminaire de Probabilit\'es XIX. Lecture Notes in Math. \textbf{1123} (1985), 130-175.

\bibitem {Bennett} G. Bennett, Schur multipliers. Duke Math J. \textbf{44} (1977), 603-639.

\bibitem {Bou} J. Bourgain, Vector-valued singular integrals and the $H_1-BMO$ duality. Probability theory and harmonic analysis. Monogr. Textbooks Pure Appl. Math. \textbf{98} (1986), 1-19.

\bibitem {BH} R. Boutonnet and C. Houdayer, Stationary characters on lattices of semisimple Lie groups. Publ. Math. Inst. Hautes Études Sci. \textbf{133} (2021), 1-46.

\bibitem {BF0} M. Bo\.zejko and G. Fendler, Herz–Schur multipliers and completely bounded multipliers of the Fourier algebra of a locally compact group. Boll. Un. Mat. Ital. A \textbf{3} (1984), 297-302.

\bibitem {BF} M. Bo\.zejko and G. Fendler, A note on certain partial sum operators. Quantum Probability. Banach Center Pubs. \textbf{73} (2006), 1-9.

\bibitem {Ca} L. Cadilhac, Noncommutative Khintchine inequalities in interpolation spaces of $L_p$-spaces. Adv. Math. \textbf{352} (2019), 265-296.

\bibitem {CCP} L. Cadilhac, J.M. Conde-Alonso and J. Parcet, Spectral multipliers in group algebras and noncommutative Calder\'on-Zygmund theory. J. Math. Pures Appl. \textbf{163} (2022), 450-472.

\bibitem {CW} L. Cadilhac and S. Wang, Noncommutative maximal ergodic inequalities for amenable groups. ArXiv: 2206.12228. 

\bibitem {CT} A.P. Calderón and A. Torchinsky, Parabolicmaximal functions associated with a distribution II. Adv. Math. \textbf{24} (1977), 101-171.

\bibitem {DCH} J. de Canni\`ere and U. Haagerup, Multipliers of the Fourier algebras of some simple Lie groups and their discrete subgroups. Amer. J. Math. \textbf{107} (1985), 455-500.

\bibitem {CR} A.I. Cano-M\'armol and \'E. Ricard, Calder\'on-Zygmund theory with noncommuting kernels via $H_1^c$. Studia Math. \textbf{277} (2024), 65-97.


\bibitem{C1} M. Caspers, A Sobolev estimate for radial $L_p$-multipliers on a class of semi-simple Lie groups. Trans. Amer. Math. Soc. \textbf{376} (2023), 
8919–8938.

\bibitem {CPPR} M. Caspers, J. Parcet, M. Perrin and \'E. Ricard, Noncommutative de Leuuw theorems. Forum Math Sigma \textbf{3} (2015), e21.

\bibitem {CPSZ} M. Caspers, D. Potapov, F. Sukochev and D. Zanin, Weak type commutator and Lipschitz estimates: resolution of the Nazarov-Peller conjecture. Amer. J. Math. \textbf{141} (2019), 593-610.

\bibitem {CaRe} M. Caspers and J. Reimann, On the best constants of Schur multipliers of second order divided difference functions. 
Math. Ann. \textbf{392} (2025), 1119-1166.

\bibitem {CS} M. Caspers and M. de la Salle, Schur and Fourier multipliers of an amenable group acting on non-commutative $L_p$-spaces. Trans. Amer. Math. Soc. \textbf{367} (2015), 6997-7013. 

\bibitem {CIOS} I. Chifan, A. Ioana, D. Osin and B. Sun, Wreath-like products of groups and their von Neumann algebras I:  $W^*$-superrigidity. Ann. of Math. \textbf{198} (2023), 1261-1303.

\bibitem {CLM} C.Y. Chuah, Z-C. Liu and T. Mei, A Marcinkiewicz multiplier theory for Schur multipliers. Anal. \& PDE \textbf{18} (2025), 1511-1530.

\bibitem {CGPT1} J.M. Conde-Alonso, A.M. Gonz\'alez-P\'erez, J. Parcet and E. Tablate, Schur multipliers in Schatten-von Neumann classes. Ann. of Math. \textbf{198} (2023), 1229-1260. 

\bibitem {CGPT2} J.M. Conde-Alonso, A.M. Gonz\'alez-P\'erez, J. Parcet and E. Tablate, A H\"ormander-Mikhlin theorem for higher rank simple Lie groups. J. Lond. Math. Soc. \textbf{111} (2025), e70137, 22 pp.

\bibitem {Co1} A. Connes, Classification of injective factors. Ann. of Math. \textbf{104} (1976), 73-115.

\bibitem {CRC} A. Connes, Classification des facteurs. In Operator algebras and applications (Kingston, Ont., 1980), Proc. Sympos. Pure Math. \textbf{38}, Amer. Math. Soc. 1982, 43-109.

\bibitem {Cor} A. C\'ordoba, The disc multiplier. Duke Math. J. \textbf{58} (1989), 21-29. 

\bibitem {CDSW} M. Cowling, B. Dorofaeff, A. Seeger and J. Wright, A family of singular oscillatory integral operators and failure of weak amenability. Duke Math. J. \textbf{127} (2005), 429-486.

\bibitem {CH} M. Cowling and U. Haagerup, Completely bounded multipliers of the Fourier algebra of a simple Lie group
of real rank one. Invent. Math. \textbf{96} (1989), 507- 549.

\bibitem {DRdF} J. Duoandikoetxea and J.L. Rubio de Francia, Estimations indépendantes de la dimension pour les transformées de Riesz. C. R. Acad. Sci. Paris \textbf{300} (1985), 193-196.

\bibitem {Fe} C. Fefferman, The multiplier problem for the ball. Ann. of Math. \textbf{94} (1971), 330-336.

\bibitem {GJP0} A. Gonz\'alez-P\'erez, M. Junge and J. Parcet, Smooth Fourier multipliers in group algebras via Sobolev dimension. Ann. Sci. Éc. Norm. Supér.  \textbf{50} (2017), 879–925.

\bibitem {GJP} A. Gonz\'alez-P\'erez, M. Junge and J. Parcet, Singular integrals in quantum Euclidean spaces. Mem. Amer. Math. Soc. \textbf{272}, 2021.

\bibitem {GPPR} A.M. Gonz\'alez-P\'erez, J. Parcet, J. P\'erez-Garc\'ia and \'E. Ricard, Riesz-Schur transforms. ArXiv: 2411.09324.

\bibitem {GPX} A.M. Gonz\' alez-P\' erez, J. Parcet and R. Xia, Noncommutative Cotlar identities for groups acting on tree-like structures. ArXiv: 2209.05298.

\bibitem {Gra} L. Grafakos, Classical Fourier Analysis. 2nd edition. Graduate Texts in Math. Springer, 2009.

\bibitem {Gr} A. Grothendieck, R\'esum\'e de la th\'eorie m\'etrique des produits tensoriels topologiques. Boll. Soc. Mat. Sao-Paulo \textbf{8} (1956), 1-79.

\bibitem {GV} R. Gundy and N. Varopoulos, Les transformations de Riesz et les int\'egrales stochastiques. C.R. Acad. Sci. Paris  \textbf{289} (1979), 13-16.

\bibitem {H} U. Haagerup, An example of a nonnuclear $\mathrm{C}^*$-algebra with the metric approximation property. Invent. Math. \textbf{50} (1979), 279-293.

\bibitem {H1} U. Haagerup, Group C$^*$-algebras without the completely bounded approximation property. J. Lie Theory 26 (2016), 861-887.

\bibitem {Ha} A. Harcharras, Fourier analysis, Schur multipliers on $S_p$ and non-commutative $\Lambda(p)$-sets. Studia Math. \textbf{137} (1999), 203-260.

\bibitem {HLX} G. Hong, X. Lai and B. Xu, Maximal singular integral operators acting on noncommutative $L_p$-spaces. Math. Ann. \textbf{386} (2023), 375-414.

\bibitem {HLW} G. Hong, B. Liao, and S. Wang, Noncommutative maximal ergodic inequalities associated with doubling conditions. Duke Math. J. \textbf{170} (2021), 205-246.

\bibitem {HLM} G. Hong, H. Liu and T. Mei, An operator-valued $T1$ theory for symmetric CZOs. J. Funct. Anal. \textbf{278} (2020), 108420, 27 pp.

\bibitem {Ho} L. H\"ormander, Estimates for translation invariant operators in $L^p$ spaces. Acta Math. \textbf{104} (1960), 93-140.

\bibitem {Hou} C. Houdayer, Noncommutative ergodic theory of higher rank lattices. ICM—International Congress of Mathematicians \textbf{4}, Sections 5–8, 3202–3223. EMS Press 2023.

\bibitem {HSV} C. Houdayer, D. Shlyakhtenko and S. Vaes, Classification of a family of non almost periodic free Araki-Woods factors. J. Eur. Math. Soc. \textbf{21} (2019), 3113-3142. 

\bibitem {Io} A. Ioana, W$^*$-superrigidity for Bernoulli actions of property (T) groups. J. Amer. Math. Soc. \textbf{24} (2011),
1175-1226.

\bibitem {IPV} A. Ioana, S. Popa and S. Vaes, A class of superrigid group von Neumann algebras. Ann. of Math. \textbf{178} (2013), 231-286.

\bibitem {JMP1} M. Junge, T. Mei and J. Parcet, Smooth Fourier multipliers on group von Neumann algebras. Geom. Funct. Anal. \textbf{24} (2014), 1913-1980. 

\bibitem {JMP2} M. Junge, T. Mei and J. Parcet, Noncommutative Riesz transforms -- Dimension free bounds and Fourier multipliers. J. Eur. Math. Soc. \textbf{20} (2018), 529-595.

\bibitem {JMPX} M. Junge, T. Mei, J. Parcet and R. Xia, Algebraic Calder\'on-Zygmund theory. Adv. Math. \textbf{376} (2021), 107443.

\bibitem {JP2} M. Junge and J. Parcet, Rosenthal's theorem for subspaces of noncommutatuve $L_p$. Duke Math. J. \textbf{141} (2008), 75-122.

\bibitem{JR} M. Junge and Z.J. Ruan, Approximation properties for noncommutative $L_p$-spaces associated with discrete groups. Duke Math. J. \textbf{117} (2003), 313-341.

\bibitem {JX} M. Junge and Q. Xu, Noncommutative maximal ergodic theorems. J. Amer. Math. Soc. \textbf{20} (2007), 385-439.

\bibitem {Ku} G. Kuhn, Convergence of Fourier series expansion related to free groups. Proc. Amer. Math. Soc. \textbf{92} (1984), 31-36.

\bibitem {KVPNAS} J. Kustermans and S. Vaes, The operator algebra approach to quantum groups. Proc. Nat. Acad. Sci. USA \textbf{97} (2000), 547-552

\bibitem {dLdlS} T. de Laat and M. de la Salle, Approximation properties for noncommutative $L_p$ of high rank lattices and nonembeddability of expanders. J. Reine Angew. Math. \textbf{737} (2018), 46-69.

\bibitem {LdlS} V. Lafforgue and M. de la Salle, Noncommutative $L_p$-spaces without the completely bounded approximation property. Duke. Math. J. \textbf{160} (2011), 71-116.

\bibitem {dL} K. de Leeuw, On $L_p$ multipliers. Ann. of Math. \textbf{81} (1965), 364-379.

\bibitem {zbMATH02681669} S. Lie, Theorie der Transformationsgruppen. Dritter (und letzten) {Abschnitt}. {Unter} {Mitwirkung} von {Fr}. {Engel} bearbeitet., Leipzig. {B}. {G}. {Teubner}. {XXVII} + 830 {S}. {{\(8^\circ\)}} (1893).

\bibitem {LP3} F. Lust-Piquard, Dimension free estimates for discrete Riesz transforms on products of abelian groups. Adv. Math. \textbf{185} (2004), 289-327.

\bibitem {Ma} V.I. Macaev, Volterra operators obtained from self-adjoint operators by perturbation. Dokl. Akad. Nauk SSSR \textbf{139} (1961),
810-813.

\bibitem {Mauc} F. Maucourant, Homogeneous asymptotic limits of Haar measures of semisimple linear groups and their lattices.
Duke Math. J. \textbf{136} (2007), 357-399.

\bibitem {Mei07} T. Mei, Operator-valued Hardy spaces. Mem. Amer. Math. Soc. \textbf{188}, 2007.

\bibitem {MR} T. Mei and \'E. Ricard, Free Hilbert transforms. Duke Math. J. \textbf{166} (2017), 2153-2182.

\bibitem {MRX} T. Mei, \'E. Ricard and Q. Xu, A H\"ormander-Mikhlin multiplier theory for free groups and amalgamated free products of von Neumann algebras. Adv. Math. \textbf{403} (2022), 108394.

\bibitem {Mey} P.A. Meyer, Transformations de Riesz pour les lois gaussiennes. S\'eminaire de Probabilit\'es XVIII. Lecture Notes in Math. \textbf{1059} (1984), 179-193.

\bibitem {Mi} S.G. Mikhlin, On the multipliers of Fourier integrals. Dokl. Akad. Nauk SSSR \textbf{109} (1956), 701-703.

\bibitem {NR} S. Neuwirth and \'E. Ricard, Transfer of Fourier multipliers into Schur multipliers and sumsets in a discrete group. Canad. J. Math. \textbf{63} (2011), 1161-1187.

\bibitem {Oz1} N. Ozawa, Solid von Neumann algebras. Acta Math. \textbf{192} (2004), 111-117.

\bibitem {OP1} N. Ozawa and S. Popa, Some prime factorization results for type II$_1$ factors. Invent. Math. \textbf{156} (2004), 223-234.

\bibitem {OP2} N. Ozawa and S. Popa, On a class of II$_1$ factors with at most one Cartan subalgebra. Ann. of Math. \textbf{172}
(2010), 713-749.

\bibitem {Pa1} J. Parcet, Pseudo-localization of singular integrals and noncommutative Calder\'on-Zygmund theory. J. Funct. Anal. \textbf{256} (2009), 509-593.


\bibitem {PRS} J. Parcet, \'E. Ricard and M. de la Salle, Fourier multipliers in $S \hskip-1pt L_n(\R)$. Duke Math. J. \textbf{171} (2022), 1235-1297.

\bibitem {PRo} J. Parcet and K. Rogers, Twisted Hilbert transforms vs Kakeya sets of directions. J. Reine Angew. Math. \textbf{710} (2016), 137-172. 

\bibitem {PST} J. Parcet, M. de la Salle and E. Tablate, The local geometry of idempotent Schur multipliers. Forum Math Pi \textbf{13} (2025), e14.

\bibitem {PT} J. Parcet and E. Tablate, Work in progress.

\bibitem {Pe} J. Peterson, $L_2$-rigidity in von Neumann algebras. Invent. Math. \textbf{175} (2009), 417-433.

\bibitem {PRiesz} G. Pisier, Riesz transforms: a simpler analytic proof of P.A. Meyer's inequality. In S\'eminaire de Probabilit\'es XXII. Lecture Notes in Math. \textbf{1321} (1988), 485-501.

\bibitem {PisAst} G. Pisier, Non-commutative vector valued $L_p$-spaces and completely $p$-summing maps. Ast\'erisque \textbf{247}. Soc. Math. France, 1998. 

\bibitem {PisSim} G. Pisier, Similarity Problems and Completely Bounded Maps. Lecture Notes in Mathematics \textbf{1618}. Springer-Verlag, 2001.

\bibitem {P2} G. Pisier, Introduction to Operator Space Theory. Cambridge University Press, 2003.

\bibitem {PisBAMS} G. Pisier, Grothendieck’s Theorem, past and present. Bull. Amer. Math. Soc. \textbf{49} (2012), 237-323. 

\bibitem {PisS} G. Pisier and D. Shlyakhtenko, Grothendieck’s theorem for operator spaces. Invent. Math. \textbf{150} (2002), 185-217.

\bibitem {PX} G. Pisier and Q. Xu, Non-commutative martingale inequalities. Comm. Math. Phys. \textbf{189} (1997),
667–698.

\bibitem {PX2} G. Pisier and Q. Xu, Non-commutative $L_p$-spaces. Handbook of the Geometry of Banach Spaces II (Eds. W.B. Johnson and J. Lindenstrauss) North-Holland (2003), 1459-1517.

\bibitem {Po1} S. Popa, On a class of type II$_1$ factors with Betti numbers invariants. Ann. of Math. \textbf{163} (2006), 809-889. 

\bibitem {Po2} S. Popa, Strong rigidity of II$_1$ factors arising from malleable actions of $w$-rigid groups I. Invent. Math. \textbf{165} (2006), 369-408.

\bibitem {Po3} S. Popa, Strong rigidity of II$_1$ factors arising from malleable actions of $w$-rigid groups II. Invent. Math. \textbf{165} (2006), 409-451.

\bibitem {PV1} S. Popa and S. Vaes, Group measure space decomposition of II$_1$ factors and W$^*$-superrigidity. Invent. Math. \textbf{182} (2010), 371-417.

\bibitem {PV2} S. Popa and S. Vaes, Unique Cartan decomposition for II$_1$ factors arising from arbitrary actions of free groups. Acta Math. \textbf{194} (2014), 237-284.

\bibitem {PS} D. Potapov and F. Sukochev, Operator-Lipschitz functions in Schatten-von Neumann classes. Acta Math. \textbf{207} (2011), 375-389.

\bibitem {RdFLP} J.L. Rubio de Francia, A Littlewood-Paley inequality for arbitrary intervals. Rev. Mat. Iberoamericana \textbf{1} (1985), 1-14.

\bibitem {RdF} J.L. Rubio de Francia, Transference principles for radial multipliers. Duke Math. J. \textbf{58} (1989), 1-19.

\bibitem {dlS} M. de la Salle, Analysis on simple Lie groups and lattices. ICM—International Congress of Mathematicians \textbf{4}, Sections 5-8, 3166–3188. EMS Press 2023.

\bibitem {Schur} J. Schur, Bemerkungen zur Theorie der beschr\"ankten Bilinearformen mit unendlich vielen Ver\"andlichen. J. Reine Angew. Math. \textbf{140} (1911), 1-28.

\bibitem {St3} E.M. Stein, Some results in harmonic analysis in $\R^n$, for $n \to \infty$. Bull. Amer. Math. Soc. \textbf{9} (1983), 71-73.

\end{thebibliography}

\enlargethispage{1cm}

\vskip-5pt

\null \hfill \textbf{Javier Parcet} \\
\null \hfill \small{\texttt{parcet@icmat.es}} \\
\null \hfill \small{Instituto de Ciencias Matem\'aticas, CSIC} \\
\null \hfill \small{Nicolás Cabrera 13-15, 28049, Madrid, Spain}
\end{document}